\documentclass[reqno, 11pt, twoside]{article}
\usepackage[a4paper,hmargin=1in,vmargin=1in]{geometry}
\usepackage{amsmath,amssymb,amsthm}
\usepackage{graphicx}
\usepackage{mathtools,bm}
\usepackage{array, booktabs}
\usepackage{etoolbox}
\usepackage{fourier}
\usepackage[lf]{ebgaramond}
\numberwithin{equation}{section}
\theoremstyle{plain}  
\newtheorem{theorem}{Theorem}[section]
\newtheorem{proposition}[theorem]{Proposition}
\newtheorem{corollary}[theorem]{Corollary}
\newtheorem{lemma}[theorem]{Lemma}

\theoremstyle{definition}
\newtheorem{definition}[theorem]{Definition}

\theoremstyle{remark}
\newtheorem{remark}[theorem]{Remark}

\newcommand{\addQEDstyle}[2]{\AtBeginEnvironment{#1}{\pushQED{\qed}\renewcommand{\qedsymbol}{#2}}\AtEndEnvironment{#1}{\popQED}}
 \addQEDstyle{example}{$\triangleleft$}
  \addQEDstyle{definition}{$\triangleleft$}
 \addQEDstyle{remark}{$\blacktriangleleft$}

\newcommand{\W}{\mathcal{W}}

\newcommand{\Sct}[1]{\mathfrak{S}_{#1}}
\newcommand{\SD}{\mathtt{SD}}
\renewcommand{\Re}{\operatorname{Re}}
\renewcommand{\Im}{\operatorname{Im}}
\newcommand{\er}{\mathrm{e}}
\newcommand{\ir}{\mathrm{i}}
\newcommand{\dr}{\mathrm{d}}
\renewcommand{\tilde}{\widetilde}
\usepackage[nobottomtitles,pagestyles]{titlesec}
\widenhead*{0.5in}{0.5in}
\renewpagestyle{plain}[\bfseries]{\footrule\setfoot{}{\thepage}{}}
\newpagestyle{mystyle}[\bfseries]{
\headrule
\sethead[\thepage][][Sloshing, Steklov and corners: Asymptotics of sloshing eigenvalues]{Michael Levitin, Leonid Parnovski, Iosif Polterovich,  and David A Sher}{}{\thepage}
}
\renewcommand\footnotemark{}
\usepackage[colorlinks,allcolors=blue]{hyperref}
\begin{document}
\pagestyle{mystyle}
\thispagestyle{plain}
\title{Sloshing, Steklov and corners:\\ 
Asymptotics of sloshing eigenvalues}
\author{Michael Levitin\hspace{-3ex}
\thanks{{\bf ML:} Department of Mathematics and Statistics,
University of Reading, Whiteknights, PO Box 220, Reading RG6 6AX, UK;
M.Levitin@reading.ac.uk; \url{http://www.michaellevitin.net}}
\and Leonid Parnovski\hspace{-3ex}
\thanks{{\bf LP:} Department of Mathematics, University College London, 
Gower Street, London WC1E 6BT, UK;
leonid@math.ucl.ac.uk}
\and Iosif Polterovich\hspace{-3ex}
\thanks{{\bf IP:} D\'e\-par\-te\-ment de math\'ematiques et de
sta\-tistique, Univer\-sit\'e de Mont\-r\'eal CP 6128 succ
Centre-Ville, Mont\-r\'eal QC  H3C 3J7, Canada;
iossif@dms.umontreal.ca; \url{http://www.dms.umontreal.ca/\~iossif}}
\and David A. Sher
\thanks{{\bf DAS:} Department of Mathematical Sciences, DePaul University, 2320 N Kenmore Ave, Chicago, IL 60614, USA;
dsher@depaul.edu}
}
\date{\small  
this version arXiv:1709.01891v5 (January 2025)\\
published version Journal d'Analyse Math\'{e}matique \textbf{146} (2022), 65--125, doi \href{https://doi.org/10.1007/s11854-021-0188-x}{10.1007/s11854-021-0188-x}}
\maketitle
\begin{abstract} In the present paper we develop an approach to obtain sharp spectral asymptotics for Steklov type problems on planar domains with corners. Our main focus is on the two-dimensional sloshing problem, which  is a mixed Steklov-Neumann boundary value problem describing small vertical oscillations of an ideal fluid in a container or in a canal with a uniform cross-section. We prove a two-term asymptotic formula for sloshing eigenvalues.  In particular, this  confirms a conjecture  posed by Fox and Kuttler in 1983.  We also obtain similar eigenvalue asymptotics for other related mixed Steklov type problems, and discuss applications to the study of Steklov spectral asymptotics on polygons.
\end{abstract}
\tableofcontents

\section{Introduction and main results}
\subsection{The sloshing problem} 
\label{section:sloshing}
Let $\Omega$ be a simply connected bounded  planar domain  with Lipschitz boundary 
such that  $\partial \Omega=S\sqcup \W$, where $S:=(A,B)$ is a line segment.  Let $0<\alpha,\beta\le\pi$ be the angles at the vertices  $A$ and $B$, respectively  (see Figure 1). 
Without loss of generality we can assume that in Cartesian coordinates $A=(0,0)$ and $B=(0,L)$, where $L>0$ is the length of $S$.

\begin{figure}[htb]
\begin{center}
\includegraphics{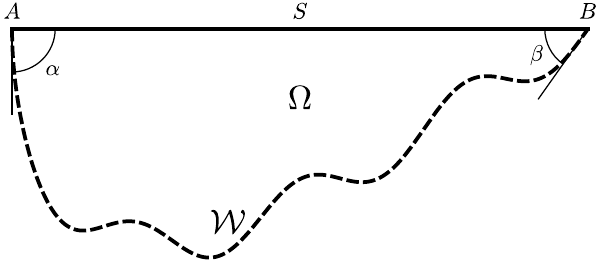}
\end{center}
\caption{Geometry of the sloshing problem\label{fig:geom}}
\end{figure}

Consider the following mixed Steklov-Neumann boundary value problem,
\begin{equation}\label{eq:sloshingproblem}
\begin{dcases} \Delta u=0 &\quad \text{in }\Omega, \\
\frac{\partial u}{\partial n}=0  &\quad \text{on } \W,\\
\frac{\partial u}{\partial n}=\lambda u &\quad \text{on }S,
\end{dcases}
\end{equation}
where $\dfrac{\partial}{\partial n}$ denotes the external normal derivative on $\partial\Omega$.
   
The eigenvalue problem \eqref{eq:sloshingproblem} is called the {\it sloshing problem}. It describes small vertical oscillations of an ideal fluid in a container or in a canal with a uniform cross-section which has the shape of the domain $\Omega$. The part $S=(A,B)$
of the boundary is called the {\it sloshing surface};  it  represents the free surface of the fluid. The part  $\mathcal{W}$  is called the {\it walls} and corresponds to the walls and the bottom of a container or a canal. The points $A$ and $B$ at the interface of the sloshing surface and the walls are called the {\it corner points}.
We also say that the walls $\W$ are {\it straight near the corners} if there exist points $A_1, B_1\in \W$ such that 
 the line segments $AA_1$ and $BB_1$ are subsets of $\W.$

\begin{figure}[htb]
\begin{center}
\includegraphics{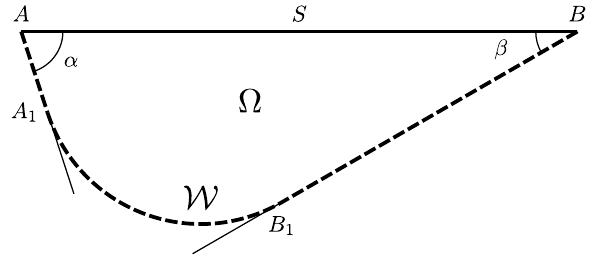}
\end{center}
\caption{Geometry of the sloshing problem with walls straight near the corners\label{fig:geom2}}
\end{figure}

It follows from general results on the Steklov type problems that the spectrum of the sloshing problem is discrete (see \cite{Ag2006, GP2017}). We denote the sloshing eigenvalues by
$$0=\lambda_1<\lambda_2 \le \lambda_3 \le \dots  \nearrow \infty,$$
where the eigenvalues are a priori counted with multiplicities.  The correspoding sloshing eigenfunctions are denoted by $u_k$, where $u_k \in C^\infty(\Omega)$, and the restrictions $u_k|_{S}$  form an orthogonal basis in $L^2(S)$.  Let us note that in two dimensions all sloshing eigenvalues are conjectured to be simple, see \cite{KKM, GP2017a}. While in full generality this conjecture remains open, in Corollary \ref{mult} we prove that it holds for all but possibly a finite number of eigenvalues.

Denote by $\mathcal{D}:L^2(S)\to L^2(S)$ the \emph{sloshing operator}, which is essentially the Dirichlet-to-Neumann operator on $S$ corresponding to Neumann boundary conditions on $\W$: given a function $f\in L^2(S)$, we have 
\begin{equation}\label{eq:DNmap}
\mathcal{D}_\Omega f=\mathcal{D}f:=\left.\frac{\partial}{\partial n}(\mathcal{H}f)\right|_S,
\end{equation}
where $\mathcal{H}f$ is the harmonic extension of $f$ to $\Omega$ with the homogeneous Neumann conditions on $\W$. Then the sloshing eigenvalues are exactly the eigenvalues of $\mathcal{D}$, and the sloshing eigenfunctions $u_k$ are harmonic extensions of the eigenfunctions of $\mathcal{D}$ to $\Omega$.

The physical meaning of the sloshing eigenvalues and eigenfunctions is as follows:  an eigenfunction  $u_k$ is the fluid velocity potential and $\sqrt{\lambda_k}$ is proportional to the frequency of the corresponding fluid oscillations. The research on the sloshing problem has a long history in hydrodynamics (see \cite[Chapter 9]{Lamb} and \cite{Greenhill}); we refer to \cite{fk} for a historical discussion.  A more recent exposition and further references could be found,  for example,  in \cite{KK, KMV, Ib, KK+}.

In this paper we investigate  the asymptotic distribution of the sloshing eigenvalues.  In fact, we develop an approach that could be applied  to study sharp spectral asymptotics of general Steklov type problems on planar domains with corners. The main difficulty is that in the presence of singularities,  the corresponding Dirichlet-to-Neumann operator is not pseudodifferential, and therefore new techniques  have to be invented, see  \cite[Section 3]{GP2017} and Subsection \ref{subs:polygons} for a discussion. Since the first version of the present  paper appeared on the arXiv, this problem has received significant  attention, both in planar and higher-dimensional cases \cite{LH, GLPS, Iv}.  Our method is based on quasimode analysis, see Subsection  \ref{subsection:outline}. A particularly challenging aspect of the argument is to show that the constructed system of quasimodes is, in an appropriate sense, {\it complete}. This is done via a rather surprising link to the theory of higher-order Sturm-Liouville problems, see Subsection \ref{stliouv}. In particular, we notice that the quasimode approximation for the sloshing problem is sensitive to the {\it arithmetic} properties of the angles at the corners. As shown in Subsection  \ref{expacc},  the quasimodes are  exponentially accurate for angles of the form $\pi/2q$, $q\in \mathbb{N}$, which together with domain monotonicity arguments allows us  to prove completeness for {\it arbitrary} angles. Further applications of our method, notably to the Steklov problem on polygons, are  presented in~\cite{LPPS}.

\subsection{Asymptotics of the sloshing eigenvalues}
As was shown in \cite{Sandgren},  if the boundary of $\Omega$ is $C^2$-regular in a neighbourhood of the corner points $A$ and $B$, then as $k\to+\infty$,
$$
\lambda_k L=\pi k + o(k).
$$
It follows from the results of \cite{Ag2006}  on Weyl's law for mixed Steklov type problems that the  $C^2$ assumption can be relaxed to $C^1$. In 1983, Fox and Kuttler used numerical evidence to conjecture that the sloshing eigenvalues of domains having both interior angles at the corner points $A$ and $B$ equal to $\alpha$,   satisfy the two-term asymptotics \cite[Conjecture 3]{fk}:
\begin{equation}
\label{FKasymp}
\lambda_k L = \pi \left(k-\frac 12\right)-\frac{\pi^2}{4  \alpha}+ o(1), 
\end{equation}  
 (note that the numeration of eigenvalues in \cite{fk} is shifted by one compared to ours).
The first main result of the present paper confirms this prediction.

\begin{theorem}\label{thm:main} 
Let $\Omega$ be a simply connected bounded Lipschitz planar domain with the sloshing surface $S=(A,B)$ of length $L$ and walls $\W$ which are $C^1$-regular near the corner points $A$ and $B$. Let  $\alpha$ and $\beta$ be the interior angles between $\W$ and $S$ at the points $A$ and $B$, resp., and assume $0<\beta \le \alpha < \pi/2 $. Then  the following asymptotic expansion holds as $k\to \infty$:
\begin{equation}
\label{eq:sloshingasymptotics}\lambda_k L=\pi\left(k-\frac 12\right)-\frac{\pi^2}{8}\left(\frac 1{\alpha}+\frac 1{\beta}\right)+r(k), \quad \text{where} \quad r(k)=o(1).
\end{equation}
If, moreover, the walls $\W$ are  straight near the corners, then 
\begin{equation}
\label{errorterm}
r(k)=O\left(k^{1-\frac{\pi}{2\alpha}}\right).
\end{equation}
\end{theorem}

In particular, for $\alpha=\beta$ we obtain \eqref{FKasymp} which proves the Fox--Kuttler conjecture for all angles $0<\alpha<\pi/2$.  Let us note that for $\alpha=\beta=\pi/2$ the asymptotics \eqref{FKasymp} have been earlier established in \cite{dav65, Ursell, dav74}. Moreover, for $\alpha=\beta=\pi/2$ it was shown that there exist further terms in the  asymptotics \eqref{FKasymp} involving  the curvature of $\W$ near the corner points. 

Before stating our next result we require the following definition.
\begin{definition} 
\label{localJohn}
A corner point $V\in\{A,B\}$ is said to \emph{satisfy a local John's condition} if  there exist a neighbourhood $\mathcal{O}_V$ of the point $V$ such that the orthogonal projection of $\W \cap\mathcal{O}_V$ onto the $x$-axis is a subset of $[A,B]$.
\end{definition}

For $\alpha=\pi/2\ge\beta$ we obtain the following
\begin{proposition}
\label{prop:pi2}
In the notation of Theorem \ref{thm:main}, let  $\alpha=\pi/2>\beta$, and assume that $A$ satisfies a local John's condition. Then
\begin{equation}
\label{pi2}
\lambda_k L = \pi \left(k-\frac 34\right)-\frac{\pi^2}{8  \beta}+ r(k), \,\,\, r(k)=o(1).
\end{equation}
The same result holds if $\alpha=\beta=\pi/2$, and additionally  $B$ satisfies a local John's condition.

If, moreover, the walls $W$ are straight near both corners, then
$$
r(k)=O\left(k^{1-\frac{\pi}{2\beta}}\right),
$$
provided $\beta<\pi/2$ and 
\[
r(k)=o\left(\er^{-k/C}\right),
\]
if $\alpha=\beta=\pi/2$.
\end{proposition}
Proposition \ref{prop:pi2} provides a solution to \cite[Conjecture~4]{fk} under the assumption that the corner point corresponding to the angle $\pi/2$ satisfies local John's condition.
\begin{remark}
In Theorem \ref{thm:main}, and everywhere further on, $C$ will denote various positive constants which depend only upon the domain $\Omega$.
\end{remark}

\begin{remark}
Definition \ref{localJohn} is a local version of John's condition which often appears  in sloshing problems, see  \cite{KKM, BKPS}. 
\end{remark}

 Theorem \ref{thm:main} and Proposition \ref{prop:pi2} yield the following
\begin{corollary}
\label{mult}
Given a sloshing problem on a domain $\Omega$ satisfying the assumptions of Theorem \ref{thm:main} or Proposition \ref{prop:pi2},  there exists $N>0$, such that 
the eigenvalues $\lambda_k$ are simple for all $k\ge N$.
\end{corollary}
This result partially confirms the conjecture about the simplicity of sloshing eigenvalues mentioned in subsection \ref{section:sloshing}.

\subsection{Eigenvalue asymptotics for a mixed Steklov-Dirichlet problem} Boundary value problems of  Steklov type with mixed boundary conditions admit several physical and probabilistic interpretations (see \cite{BK, BKPS}).  In particular, the sloshing problem \eqref{eq:sloshingproblem} could be also used to model
 the stationary heat distribution in $\Omega$ such that the walls $\W$ are perfectly insulated and the heat flux through $S$ is proportional to the temperature.
If instead the walls  $\W$ are  kept under zero temperature, one obtains the following mixed Steklov-Dirichlet problem:
\begin{equation}\label{eq:SteklovDirichlet}
\begin{dcases} 
\Delta u=0 & \quad\text{in }\Omega, \\
u=0  &\quad \text{on } \W,\\ 
\frac{\partial u}{\partial n}=\lambda^D u &\quad \text{on }S,
\end{dcases}
\end{equation}
Let $0<\lambda_1^D \le \lambda_2^D \le \dots \nearrow \infty$ be the eigenvalues of the problem \eqref{eq:SteklovDirichlet}. 
Similarly to Theorem \ref{thm:main} we obtain
\begin{theorem}\label{thm:Dirichlet} 
Assume that the domain $\Omega$ and its boundary $\partial \Omega=S \sqcup \W$  satisfy the assumptions of  Theorem \ref{thm:main}.  Let  $\alpha$ and $\beta$ be the interior angles between $\W$ and $S$ at the points $A$ and $B$, resp., and assume $0<\beta \le \alpha < \pi/2 $. Then  the following asymptotic expansion holds as $k\to \infty$: 
\begin{equation}\label{eq:sloshingasymptoticsdirichlet}\lambda_k^D L=\pi\left(k-\frac 12\right)+\frac{\pi^2}{8}\left(\frac 1{\alpha}+\frac 1{\beta}\right)+r^D(k),
\quad \operatorname{where} \,\,\,r^D(k)=o(1).
\end{equation}
If, moreover, the walls $\W$ are straight near the corner points $A$ and $B$, then
\begin{equation}
\label{errortermdir}
r^D(k)= O\left(k^{1-\frac{\pi}{\alpha}}\right).
\end{equation}
\end{theorem}

We also have the following analogue of Proposition \ref{prop:pi2}.
\begin{proposition}
\label{prop:pi22}
In the notation of Theorem \ref{thm:Dirichlet}, let  $\alpha=\pi/2 > \beta$, and assume that $A$ satisfies a local John's condition. Then
\begin{equation}
\label{pi22}
\lambda_k^D L = \pi \left(k-\frac 14\right)+\frac{\pi^2}{8  \beta}+ r(k), \,\,\, r(k)=o(1).
\end{equation}
The same result holds if $\alpha=\beta=\pi/2$ and additionally $B$ satisfies a local John's condition.

If, moreover, the walls $W$ are straight near both corners, then
$$
r(k)=O\left(k^{1-\frac{\pi}{\beta}}\right).
$$
provided $\beta<\pi/2$, and 
$$
r(k)=o\left(\er^{-k/C}\right)
$$
for some $C>0$ if $\alpha=\beta=\pi/2$.
\end{proposition}

This result will be used for our subsequent analysis of polygonal domains, see Subsection \ref{subs:polygons}.

The analogue of Corollary \ref{mult} clearly holds in the Steklov-Dirichlet case as well.
\begin{remark}
\label{belief}
 We believe that 
asymptotic formulae \eqref{eq:sloshingasymptotics} and \eqref{eq:sloshingasymptoticsdirichlet} in fact hold for any angles $\alpha,\beta  \le \pi$. 
Note also  that if the walls are straight near the corners, our method yields a slightly better  remainder estimate for the Steklov-Dirichlet problem compared to the Steklov-Neumann one.  We refer to the proof of Theorem \ref{wedgemodels} for details.
\end{remark}

\subsection{Outline of the approach}
\label{subsection:outline}
Let us sketch the main ideas of the proof of Theorem \ref{thm:main}; modifications needed to obtain Theorem \ref{thm:Dirichlet} are quite minor. First, we observe that using domain monotonicity for sloshing eigenvalues (see \cite{BKPS}), one can deduce the general asymptotic expansion \eqref{eq:sloshingasymptotics}  from the two-term asymptotics with the remainder \eqref{errorterm} for domains with straight walls near the corners.
Assuming that the walls are straight near the corner points, we  construct the \emph{quasimodes}, i.e. the  approximate eigenfunctions of the problem \eqref{eq:sloshingproblem}. This is done by transplanting certain model solutions of the mixed Steklov-Neumann  problem in an infinite angle (cf. \cite{DW} where a similar idea has been implemented at a physical level of rigour). These solutions are in fact of independent interest: they were used to describe ``waves on a sloping beach" (see \cite{Hanson, lewy, Sto47, peters, Ehr87}).  The approximate eigenvalues given by the first two terms on the right hand side of \eqref{eq:sloshingasymptotics} are then found from a matching condition between the two model solutions transplanted to the corners $A$ and $B$, respectively. Given that the model solutions decay rapidly away from the sloshing surface, it follows  that the shape of the walls away from the corners does not matter for our approximations, so in fact the domain $\Omega$ can be viewed as a triangle with angles $\alpha$ and $\beta$ at the sloshing surface $S$. From the standard quasimode analysis it follows that there exist a sequence of eigenvalues of the sloshing satisfying the asymptotics \eqref{eq:sloshingasymptotics}. A major remaining challenge is to show that this sequence is {\it asymptotically complete} (i.e. that there are no other sloshing  eigenvalues that have not been accounted for, see subsection \ref{compabst} for a formal definition),  and that the enumeration of eigenvalues given by \eqref{eq:sloshingasymptotics} is correct.
This can not be achieved by simple arguments. While the set of quasimodes is a  perturbation of a set of exponentials, even if one can prove the completeness of this latter set, the perturbation is too large to apply the standard Bary-Krein result (Lemma \ref{lem:barykrein}) to deduce the claimed asymptotic completeness. 

Our quasimode construction for arbitrary angles $\alpha, \beta < \pi/2$ is based on the model solutions obtained by Peters \cite{peters}. These solutions are given in terms of some complex integrals, and while their asymptotic representations allow us to construct the quasimodes, they are not accurate enough to prove completeness. However, it turns out that for angles $\alpha=\beta=\pi/2q$,  $q\in \mathbb{N}$,  model solutions can be written down explicitly as linear combinations of certain complex exponentials. Moreover, it turns out that for such angles the model solutions can be used to approximate the eigenfunctions of a Sturm-Liouville type problem of order $2q$ with Neumann boundary conditions  (see Theorem \ref{ODE}), for which the completeness follows from the general theory of linear ODEs. The enumeration of the sloshing eigenvalues in \eqref{eq:sloshingasymptotics} may then be verified  by developing the approach outlined  in \cite[Chapter 2]{naimark}. Another important property used here is a peculiar duality between the Dirichlet and Neumann eigenvalues of the Sturm-Liouville problem, see Proposition \ref{prop:weirdODEfact}.

Once we have proven completeness and established  the enumeration of the sloshing eigenvalues in the case $\alpha=\beta=\pi/2q$, we use once again domain monotonicity together with  continuous perturbation arguments to complete the proof of Theorem \ref{thm:main} for arbitrary angles.
\subsection{An application to higher order Sturm-Liouville problems}
\label{stliouv}
Let us elaborate on  the link between the sloshing problem with angles  $\alpha=\beta=\pi/2q$ and the higher order Sturm-Liouville type ODEs mentioned in
the previous section. 
For a given $q\in \mathbb{N}$, consider an eigenvalue problem on an interval $(A,B) \subset \mathbb{R}$ of length $L$:
\begin{equation}
\label{eq:ode}
(-1)^q U^{(2q)}(x) = \Lambda^{2q} U(x),
\end{equation}
with either Dirichlet
\begin{equation}
\label{bc:Dirichlet}
U^{(m)}(A)=U^{(m)}(B)=0,\,\,  m=0,1,\dots, q-1, 
\end{equation}
or Neumann
\begin{equation}
\label{bc:Neumann}
U^{(m)}(A)=U^{(m)}(B)=0, \,\, m=q,q+1,\dots,2q-1, 
\end{equation}
boundary conditions.
For  $q=1$ we obtain the classical Sturm-Liouville equation describing vibrations of either a fixed or a free string. The case $q=2$ yields the vibrating beam equation, also with either fixed or free ends. It follows from general elliptic theory that the spectrum of the boundary value problems \eqref{bc:Dirichlet} or  \eqref{bc:Neumann} for the equation \eqref{eq:ode} is discrete. It is easy to check that all Dirichlet eigenvalues are positive, while the Neumann spectrum contains an eigenvalue zero of multiplicity $q$; the corrresponding eigenspace is generated by the functions $1, x,\dots, x^{q-1}$.   Let
$$\underbrace{0=\dots...=0}_{q\,\,  \operatorname{times}}<(\Lambda_{q+1})^{2q}\le (\Lambda_{q+2})^{2q} \le \dots \nearrow \infty$$ be the spectrum of the Neumann problem \eqref{bc:Neumann}. Then, as shown in Proposition \ref{prop:weirdODEfact}, $\Lambda_k^D=\Lambda_{k+q}$, $k=1,2,\dots$, where $(\Lambda_k^D)^{2q}$ are the eigenvalues of the Dirichlet problem \eqref{bc:Dirichlet}. 
We have the following result:
\begin{theorem}
\label{ODE}
 For any $q\in \mathbb{N}$, the following asymptotic formula holds for the eigenvalues $\Lambda_k$:
\begin{equation}
\label{precise}
\Lambda_k L= \pi \left(k-\frac{1}{2}\right) -\frac{\pi q}{2} + O\left(\er^{-k/C}\right)
\end{equation}
where $C>0$ is some positive constant. Moreover, let $\lambda_k$, $k=1,2,\dots$, be the eigenvalues of a sloshing problem \eqref{eq:sloshingproblem}
with  straight walls  near the corners making equal angles $\pi/2q$ at both corner points with the sloshing surface $S$ of length $L$. Then
\begin{equation}
\label{link}
\lambda_k  = \Lambda_k  + O\left(\er^{-k/C}\right),
\end{equation}
and the eigenfunctions $u_k$ of the sloshing problem decrease exponentially away from the sloshing boundary $S$. 
\end{theorem}
Spectral asymptotics of Sturm-Liouville type problems of arbitrary order for general self-adjoint boundary conditions have been studied earlier, see  \cite[Chapter II, section 9]{naimark} and references therein. However, for the special case of problem \eqref{eq:ode} with Dirichlet or Neumann boundary conditions, formula \eqref{precise} gives a more precise result. First, \eqref{precise} gives the asymptotics for each $\Lambda_k$, while earlier results yield only an asymptotic form of the  eigenvalues without specifying the correct numbering. Second, we get an exponential error estimate, which is an improvement upon a $O(1/k)$ that was known previously. Finally, and maybe most importantly, \eqref{link} provides a physical meaning to the Sturm-Liouville problem \eqref{eq:ode} for arbitrary order $q\ge 1$.
We also note that for $q=2$, 
the Sturm-Liouville eigenfunctions $U_k(x)$, $k=1,2,\dots$,  are precisely the traces of the sloshing eigenfunctions $\left.u_k\right|_S$ on an isosceles right triangle  $\Omega=T$ such that the sloshing surface $S=(A,B)$ is its hypotenuse. This  fact was already known to H. Lamb back in the nineteenth century (see \cite{Lamb}), and Theorem \ref{ODE} extends the connection  between higher order Sturm-Liouville problems and the sloshing problem to $q\ge 2$. Let us also note  that in a different context related to the study of photonic crystals, a connection between higher order ODEs and Steklov-type problems on domains with corners has been explored in \cite{Kuku}.

\subsection{Towards sharp asymptotics for Steklov eigenvalues on polygons}\label{subs:polygons} As was discussed in \cite[Section 3]{GP2017}, precise asymptotics of Steklov eigenvalues on polygons  and on smooth planar domains are quite different. Moreover, the powerful pseudodifferential methods used in the smooth case can not be applied for polygons, since the Dirichlet-to-Neumann operator on the boundary of a non-smooth domain is not pseudodifferential. It turns out that the methods of the present paper may be developed in order to investigate the Steklov spectral asymptotics on polygonal domains. This is the subject of  the subsequent paper \cite{LPPS}, see also \cite{KLPPS21}. 
While establishing sharp eigenvalue asymptotics for polygons requires a lot of further work, a sharp remainder estimate in Weyl's law follows immediately from Theorems \ref{thm:main} and \ref{thm:Dirichlet}.
\begin{corollary}
Let $N_\mathcal{P}(\lambda)=\#\{\lambda_k < \lambda\}$ be the counting function of Steklov eigenvalues $\lambda_k$ on a convex polygon $\mathcal P$ of perimeter $L$. Then
\begin{equation}
\label{weylpolygon}
N_\mathcal{P}(\lambda)=\frac{L}{\pi} \lambda + O(1).
\end{equation}
\end{corollary}
\begin{remark} The asymptotic formula \eqref{weylpolygon} improves upon the previously known error bound $o(\lambda)$  (see \cite{Sandgren, Ag2006}). Note that the $O(1)$ estimate for the error term in Weyl's law is {\it optimal}, since the counting function is a step-function.
\end{remark}
\begin{proof} 
Given a convex $n$-gon $\mathcal{P}$, take an arbitrary point $O \in \mathcal{P}$.  It can be connected  with the vertices of $\mathcal{P}$ by $n$ smooth curves having only point $O$ in common in such a way that at each vertex, the angles between the sides of the polygon and the corresponding curve are smaller than $\pi/2$. This is clearly possible since all the angles of a convex polygon are less than $\pi$. Let $\mathcal{L}$ be the union of those curves. 
\begin{figure}[htb]
\begin{center}
\includegraphics{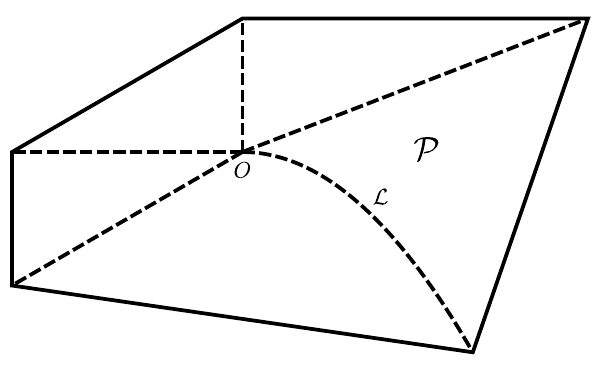}
\end{center}
\caption{A polygon with auxiliary curves  $\mathcal{L}$ \label{fig:poly}}
\end{figure}
Consider two auxiliary spectral problems: in the first one, we impose Dirichlet conditions on 
$\mathcal{L}$ and keep the Steklov condition on the boundary of the polygon, and in the second one we impose Neumann conditions on $\mathcal{L}$  and keep the Steklov condition on the boundary. Let $N_1(\lambda)$ and $N_2(\lambda)$ be, respectively, the counting functions of the first and the second problem. By Dirichlet--Neumann bracketing (which works for the sloshing problems via the variational principle in the same way it does for the Laplacian) we get
$$
N_1(\lambda) \le N_\mathcal{P}(\lambda) \le N_2(\lambda), \,\,\, \lambda>0.
$$
The spectrum of the second auxiliary problem can be represented as the union of spectra of  $n$ sloshing problems, while simultaneously the spectrum of the first auxiliary problem can be represented as the union of spectra of $n$ corresponding Steklov-Dirichlet problems. Applying Theorems \ref{thm:main} and \ref{thm:Dirichlet} to those problems, and transforming the eigenvalue asymptotics into the asymptotics of counting functions, we obtain
$N_2(\lambda)-N_1(\lambda)=O(1)$, which implies \eqref{weylpolygon}. This completes the proof of the corollary.
\end{proof}
Let us conclude this subsection by a result in the spirit of Theorems \ref{thm:main} and \ref{thm:Dirichlet} that will be used in \cite{LPPS} in the proof of the sharp asymptotics of Steklov eigenvalues on polygons.
\begin{proposition} 
\label{thm:mixed}
Let $\Omega=\bigtriangleup ABZ$ be a triangle with angles $\alpha, \beta \le \pi/2$ at the vertices $A$ and $B$, respectively. Consider a mixed Steklov-Dirichlet-Neumann problem on this triangle with the Steklov condition imposed on $AB$, Dirichlet condition imposed on $AZ$ and Neumann condition imposed on $BZ$. Assume that the side $AB$ has length $L$. Then the eigenvalues $\lambda_k$ of the mixed Steklov-Dirichlet-Neumann problem on $\bigtriangleup ABZ$ satisfy the asymptotics:
\begin{equation}
\label{mixedasymp}
\lambda_k L=\pi\left(k-\frac 12\right)+\frac{\pi^2}{8}\left(\frac 1{\alpha}-\frac 1{\beta}\right)+O\left(k^{1-\frac{\pi}{\max(\alpha,2\beta)}}\right).
\end{equation}
\end{proposition}

\begin{remark} 
Note that the Dirichlet condition near the vertex $A$ yields a contribution $\dfrac{\pi^2}{8\alpha}$ (with a positive sign) to the eigenvalue asymptotics, while
the Neumann condition near $B$ contributes $-\dfrac{\pi^2}{8\beta}$ (with a negative sign). This is in good agreement with the intuition provided by Theorems \ref{thm:main} and \ref{thm:Dirichlet}.  
\end{remark}

\begin{remark}\label{rem:believe}
In fact, we believe that a stronger statement than the one proposed in Remark \ref{belief} is true:  formula \eqref{mixedasymp} holds for any mixed Steklov-Dirichlet-Neumann problem on a domain with a curved Steklov part $AB$ of length $L$, curved walls $\W$, and angles $\alpha,\beta<\pi$ at $A$ and $B$, such that the Dirichlet condition is imposed near $A$ and Neumann condition is imposed near $B$.
\end{remark}

\subsection{Plan of the paper}\label{plan} In Section \ref{sec:triangular}  we use the Peters solutions of the sloping beach problem \cite{peters} to construct quasimodes for the sloshing and Steklov-Dirichlet problems on triangular domains. In Section \ref{sec:expquasimodes} we consider the case of angles of the form $\pi/2q$ for some positive integer $q$.
In Section \ref{sec:completeness} we first prove the completeness of this system of exponentially accurate quasimodes using a connection to higher order Sturm-Liouville eigenvalue problems. In particular, we prove Theorem \ref{ODE}. After that, we apply domain monotonicity arguments in order to prove Theorems \ref{thm:main} and \ref{thm:Dirichlet}, as well as Propositions \ref{prop:pi2}, \ref{prop:pi22} and \ref{thm:mixed} for triangular domains. In section \ref{sec:curvilinear} we extend these results to domains with curvilinear walls: first, for domains with the walls  that are straight near the corners, and then to domains with general curvilinear walls.  In Appendix \ref{proofpeters} we prove Theorem \ref{wedgemodels}, which is essentially based on the ideas of \cite{peters}. In Appendix \ref{naimarkproof} we prove an auxiliary Proposition \ref{naimarklemma} which is needed to prove Theorem \ref{ODE}. This section draws heavily on the results of  \cite[Chapter 2]{naimark}. Some numerical examples are presented in Appendix \ref{app:examples}. 

\subsection*{Acknowledgments} The authors are grateful to Lev Buhovski for suggesting the approach used in section \ref{subs:curvilin}, to Rinat Kashaev for proposing an alternative route to prove Theorem \ref{thm:derivsofhlsols},  as well as to Alexandre Girouard and Yakar Kannai for numerous useful discussions on this project.  We are also grateful to Marcello Malagutti for his comments on the proof of Lemma A.4, see a footnote on page \pageref{page:foot}.
The research of L.P. was supported by by EPSRC grant EP/J016829/1.
The research of I.P. was supported by NSERC, FRQNT, Canada Research Chairs program, as well as the Weston Visiting Professorship program  at  the Weizmann Institute of Science,
where part of this work has been accomplished.
The research of D.S. was supported in part by NSF EMSW21-RTG 1045119 and in part by a Faculty Summer Research Grant from DePaul University.

\section{Construction of quasimodes for triangular domains}\label{sec:triangular}
\subsection{The sloping beach problem} Let $(x,y)$ be Cartesian coordinates in $\mathbb{R}^2$, let $z=x+\ir y\in\mathbb{C}$, and let $(\rho,\theta)$ denote the polar coordinates $z=\rho \er^{\ir \theta}$. Let $\Sct{\alpha}$ be the planar sector $-\alpha\leq\theta\leq 0$. Consider the following mixed Robin-Neumann problem
\begin{equation}\label{modelproblem:neumann}
\begin{dcases}\Delta\phi=0&\quad\text{in }\Sct{\alpha},\\ 
\frac{\partial\phi}{\partial y}=\phi&\quad\text{on }\Sct{\alpha}\cap\{\theta=0\},\\
\frac{\partial\phi}{\partial n}=0&\quad\text{on }\Sct{\alpha}\cap\{\theta=-\alpha\},\\
\end{dcases}
\end{equation}
and the mixed Robin-Dirichlet problem

\begin{equation}\label{modelproblem:dirichlet}
\begin{dcases}\Delta\phi=0&\quad\text{in }\Sct{\alpha},\\ 
\frac{\partial\phi}{\partial y}=\phi&\quad\text{on }\Sct{\alpha}\cap\{\theta=0\},\\
\phi=0&\quad\text{on }\Sct{\alpha}\cap\{\theta=-\alpha\}.\\
\end{dcases}
\end{equation}
Note that there is no spectral parameter in these problems (hence the boundary conditions are called Robin rather than Steklov). 
Our aim is to exhibit solutions of \eqref{modelproblem:neumann} and \eqref{modelproblem:dirichlet} decaying away from the horizontal part of the boundary with the property that $\phi(x,0)\to\cos(x-\chi)$ for some fixed $\chi$ as $x\to\infty$. 

This problem is known as the {\it sloping beach problem}, and has a long and storied history in hydronamics, see \cite{lewy} references therein. It turns out that the form of the solutions depends in a delicate way on the \emph{arithmetic} properties of the angle $\alpha$; we will discuss this in more detail later on. Let 
\begin{equation}\label{eq:muchi}
\mu_\alpha=\frac{\pi}{2\alpha},\quad \chi_{\alpha,N}=\frac{\pi}{4}(1-\mu_\alpha),\quad \chi_{\alpha,D}=\frac{\pi}{4}(1+\mu_\alpha).
\end{equation}
The following key result was essentially established by Peters \cite{peters}:

\begin{theorem}\label{wedgemodels} For any $0<\alpha<\pi/2$,  there exist solutions $\phi_{\alpha,N}(x,y)$ and $\phi_{\alpha,D}(x,y)$ of \eqref{modelproblem:neumann} and \eqref{modelproblem:dirichlet},  respectively, and a constant $C>0$ such that:
\begin{itemize}
\item $|\phi_{\alpha,N}(x,y)|$ and $|\phi_{\alpha,D}(x,y)|$ are bounded on the closed sector $\overline{\Sct{\alpha}}$;
\item 
\begin{equation}\label{eq:phialphaN}
\phi_{\alpha,N}(x,y)=\er^{y}\cos(x-\chi_{\alpha,N})+R_{\alpha,N}(x,y),
\end{equation} 
where
\[
\left|R_{\alpha,N}(x,y)\right|\leq C\rho^{-\mu},\quad \left|\nabla_{(x,y)}R_{\alpha,N}(x,y)\right|\leq C\rho^{-\mu-1};
\]
\item 
\begin{equation}\label{eq:phialphaD}
\phi_{\alpha,D}(x,y)=\er^{y}\cos(x-\chi_{\alpha,D})+R_{\alpha, D}(x,y),
\end{equation}
where
\[
\left|R_{\alpha,D}(x,y)\right|\leq C\rho^{-2\mu},\quad \left|\nabla_{(x,y)}R_{\alpha,D}(x,y)\right|\leq C\rho^{-2\mu-1};
\]
\item if $\mathcal P$ is any differential operator of order $k$ with constant coefficient 1, then as $\rho\to 0$,
\[\left|\mathcal P\phi_{\alpha,N}(x,y)\right|=o(\rho^{-k}),\quad \left|\mathcal P\phi_{\alpha,D}(x,y)\right|=O(\rho^{\mu-k}),\]
and for all $\rho$, most importantly $\rho\geq 1$,
\[\left|\mathcal P R_{\alpha,N}(x,y)\right|\leq C_{k}\rho^{-\mu-k},\quad \left|\mathcal P R_{\alpha,D}(x,y)\right|\leq C_{k}\rho^{-2\mu-k}.\]
\end{itemize}

For $\alpha=\pi/2$  we have
\begin{equation}
\label{pi2solutions}
\phi_{\frac{\pi}{2},N}=\er^y\cos x, \,\,  \phi_{\frac{\pi}{2},D}=\er^y \sin x.
\end{equation}
\end{theorem}

\begin{remark} Note that $R_{\alpha,N}(x,y)$ and $R_{\alpha,D}(x,y)$ are harmonic. Note also that our further analysis (specifically, the asympotic completeness of quasimodes) will imply that the solutions $\phi_{\alpha,N}(x,y)$ and $\phi_{\alpha,D}(x,y)$ are unique in the sense that we cannot replace the particular constants $\chi_{\alpha,N}$ and $\chi_{\alpha,D}$ from \eqref{eq:muchi} by any other value.
\end{remark}

The construction of the solutions for both the Robin-Neumann problem and the Robin-Dirichlet problem is due to Peters \cite{peters}. The approach is based on complex analysis, specifically the Wiener-Hopf method. The solution is written down explicitly as a complex integral, which allows us to analyse the asymptotics. We have reproduced this construction, taking special care with the remainder estimates that were not worked out in \cite{peters}. The proof of Theorem \ref{wedgemodels} is quite technical and is deferred until Appendix \ref{proofpeters}.

\subsection{Quasimode analysis}\label{subsection:step2} In what follows, we focus on the proof of Theorem \ref{thm:main}; minor modifications required for Theorem \ref{thm:Dirichlet}
will be discussed later. As such we suppress all $N$ subscripts.

As was indicated in subsection \ref{plan}, we split the proof of Theorem \ref{thm:main} into several steps. 
Our first  objective  is to prove the asymptotic expansion \eqref{eq:sloshingasymptotics} with the remainder estimate \eqref{errorterm} for  triangular domains.

\smallskip

\noindent{\bf Convention:}  {\it From now on and until subsection \ref{subsection:monotonicity} we assume that $\Omega$ is a triangle $T=\bigtriangleup ABZ$.}

\smallskip

The key starting idea is to glue together two scaled Peters solutions  $\pm \phi_{\alpha}(\sigma x,\sigma y)$, one at each corner $A$ and $B$, to construct quasimodes. These Peters solutions must match, meaning that  the phases of the trigonometric functions in the asymptotics of both solutions must agree. Denoting those phases by $\chi_{\alpha}=\chi_{\alpha,N}$ and $\chi_{\beta}=\chi_{\beta,N}$,  see \eqref{eq:muchi}, we require
\[
\cos(\sigma x-\chi_{\alpha})=\pm\cos((L-x)\sigma-\chi_{\beta}).
\]
Solving this equation immediately  gives $\sigma=\sigma_k$ for some integer $k$, where 
\begin{equation}\label{eq:quantizationcondition}
\sigma_kL=\pi\left(k-\frac 12\right)-\frac{\pi^2}{8}\left(\frac{1}{\alpha}+\frac{1}{\beta}\right).
\end{equation}
This could be viewed as  the {\it quantization condition} resulting in the asymptotics \eqref{eq:sloshingasymptotics}.

Let us  now make this precise. Let the plane wave $p_{\sigma}(x,y)=\er^{\sigma y}\cos(\sigma x-\chi_{\alpha})$, where $\sigma$ is determined by the quantization condition
\eqref{eq:quantizationcondition}. 
Letting $z=(x,y)$, we set
\[
R_{A,\sigma}(z):=\phi_{\alpha}(\sigma z)-p_{\sigma}(z);
\]
that is, $R_{A,\sigma}(z)$ is the difference between the scaled Peters solution in the sector of angle $\alpha$ with the vertex at $A$, and the scaled trigonometric function. Similarly, let 
\[
R_{B,\sigma}(z):=\phi_{\beta}(\sigma (L-\overline{z}))-p_{\sigma}(L-\overline{z})
\] 
be the difference between the scaled Peters solution  in the sector of angle $\beta$ with the vertex at $B$ and a scaled trigonometric function. Note that $R_{A,\sigma}(z)=R_{\alpha}(z \sigma)$ and similarly for $R_{B,\sigma}(z)$. By the respective remainder estimates in Theorem \ref{wedgemodels},
\begin{align*}
|R_{A,\sigma}(z)|&\leq C\sigma^{-\mu_\alpha}|z-A|^{-\mu_\alpha};& 
|R_{B,\sigma}(z)|&\leq C\sigma^{-\mu_\beta}|L-\overline{z}-B|^{-\mu_\beta};\\
|\nabla R_{A,\sigma}(z)|&\leq C\sigma^{-\mu_\alpha}|z-A|^{-\mu_\alpha-1};& 
|\nabla R_{B,\sigma}(z)|&\leq C\sigma^{-\mu_\beta}|L-\overline{z}-B|^{-\mu_\beta-1}.
\end{align*}
For $\sigma=\sigma_k$ satisfying the quantization condition \eqref{eq:quantizationcondition}, let us define a function $v_\sigma'$ (note that $v_\sigma'$  is not a derivative of $v_\sigma$ but a new function, and  we will follow this convention in the sequel):
\begin{equation}\label{eq:uncorrectedquasimodes}
v'_{\sigma}(z):=p_{\sigma}(z)+R_{A,\sigma}(z)+R_{B,\sigma}(z)
\end{equation}
(indeed, this definition is meaningful only when $\sigma$ satisfies \eqref{eq:quantizationcondition}). This is our first attempt at a quasimode. The problem with this function is that, while it is harmonic, it does not satisfy the Neumann condition on $\W$. However, the error is small and we can correct it.

Indeed, in a neighbourhood of $A$, the function $p_{\sigma}+R_{A,\sigma}$ is the Peters solution and hence satisfies the Neumann condition on $AZ$, and $\nabla R_{B,\sigma}$ is of order $O\left(\sigma^{-\mu_\beta}\right)$, so the normal derivative of $v'_{\sigma}$ is $O\left(\sigma^{-\mu_\beta}\right)$. A similar analysis holds in a neighbourhood of $B$. Away from both $A$ and $B$, all three terms have gradients of magnitude $O\left(\sigma^{-\mu_\alpha}\right)$ (as $\alpha\ge\beta$). Therefore
\begin{equation}\label{eq:gradcorrect}\left|\frac{\partial v'_{\sigma}}{\partial n}\right|\leq C\sigma^{-\mu_\alpha}\quad\text{ on }\W.
\end{equation}

In order to correct this ``Neumann defect", consider a function $\eta_{\sigma}$ defined as a solution of the following system:
\begin{equation}\label{eq:neumanndefect}
\begin{dcases}\Delta\eta_{\sigma}=0 &\quad\text{in }\Omega;\\ 
\frac{\partial}{\partial n}\eta_{\sigma}=\frac{\partial}{\partial n}v'_{\sigma} &\quad\text{on }\W;\\ 
\frac{\partial}{\partial n}\eta_{\sigma}=-\kappa_{\sigma}\psi &\quad\text{on }S,
\end{dcases}
\end{equation}
where $\psi\in C^{\infty}(S)$ is a fixed function, supported away from the vertices, with $\int_S \psi=1$, and where
\[
\kappa_{\sigma}=\int_{\W}\frac{\partial v'_{\sigma}}{\partial n}.
\]
Note that the integral of the Neumann data over $\partial\Omega$ in \eqref{eq:neumanndefect} is zero; thus a solution $\eta_{\sigma}$ to \eqref{eq:neumanndefect} exists and is uniquely defined up to an additive constant. It is well-known (see, for instance, \cite{Rondi}) that the Neumann-to-Dirichlet map is a bounded operator on $L^2_*(\partial\Omega)$, the space of mean-zero $L^2$ functions on the boundary. Therefore, $\|\eta_{\sigma}\|_{L^2(\partial\Omega)}\leq C\sigma^{-\mu_\alpha}$, and hence
\begin{equation}\label{etaestimate}
\|\eta_{\sigma}\|_{L^2(S)}\leq C\sigma^{-\mu_\alpha}.
\end{equation}
With this auxiliary function, we define a corrected quasimode:
\[v_{\sigma}(z):=v'_{\sigma}(z)-\eta_{\sigma}(z)=p_{\sigma}(z)+R_{A,\sigma}(z)+R_{B,\sigma}(z)-\eta_{\sigma}(z).\]
Observe that $v_{\sigma}$ is harmonic and satisfies the homogeneous Neumann boundary condition on $\W$.

The key property of our new quasimodes is the following
\begin{lemma}\label{lem:keyquasi} With the notation as above, there exists a constant $C$ such that
\begin{equation}\label{keyquasiest}
\left\|\mathcal{D}v_{\sigma}-\sigma v_{\sigma}\right\|_{L^2(S)}\leq C\sigma^{1-\mu_\alpha},
\end{equation}
where $\mathcal{D}$ is the sloshing operator \eqref{eq:DNmap}.
\end{lemma}

\begin{remark} Note that $\mu_\alpha=\frac{\pi}{2\alpha}>1$ for $\alpha<\pi/2$, and therefore we get quasimodes of order $O\left(\sigma^{-\delta}\right)$ for  $\delta:=1-\mu_\alpha>0$. For $\alpha\geq\pi/2$ one would need to modify our approach.
\end{remark}

\begin{proof}[Proof of Lemma \ref{lem:keyquasi}] Since $v_{\sigma}$ is harmonic and satisfies the Neumann condition on $\W$, we have
\[
\mathcal{D}\left(v_{\sigma}|_S\right)=\left.\frac{\partial v_{\sigma}}{\partial n}\right|_S.
\]
Consider first a region away from the vertex $A$. In this region,
\[
\left.\frac{\partial(p_{\sigma}+R_{B,\sigma})}{\partial n}\right|_S-\sigma(p_{\sigma}+\left.R_{B,\sigma})\right|_S=0
\]
by Theorem \ref{wedgemodels}. Moreover, due to the bounds on $R_{A,\sigma}$ in Theorem \ref{wedgemodels},
\[
\left|\frac{\partial R_{A,\sigma}}{\partial n}-\sigma R_{A,\sigma}\right|\leq C\sigma^{1-\mu_\alpha}\quad\text{on }S
\]
pointwise, and hence the same estimate holds  in $L^2(S)$. Finally, by construction of $\eta_{\sigma}$, we know 
\[
\left.\frac{\partial\eta_{\sigma}}{\partial n}\right|_S=-c_{\sigma}\psi,\quad |c_{\sigma}|\leq C\sigma^{-\mu_\alpha};
\]
combining this with \eqref{etaestimate} yields
\[
\left\|\frac{\partial\eta_{\sigma}}{\partial n}-\sigma\eta_{\sigma}\right\|_{L^2(S)}\leq C\sigma^{1-\mu_\alpha}.
\]
Putting everything together using the definition of $v_{\sigma}$ gives us the required  estimate  away from $A$. A similar analysis shows that the same estimate holds away from $B$, completing the proof.
\end{proof}

\begin{remark} For the Dirichlet or mixed problems on triangles, it is also possible to construct $\eta_{\sigma}(z)$ harmonic and satisfying \eqref{etaestimate} such that $v_{\sigma}(z)$ satisfies the appropriate homogeneous boundary conditions on $\W:=\W_1\cup\W_2$. In this case,  $\eta_{\sigma}$ is a solution of the following problem:
\begin{equation}\label{eq:mixeddefect}
\begin{dcases}
\Delta\eta_{\sigma}=0 &\quad\text{in }\Omega;\\ 
\frac{\partial}{\partial n}\eta_{\sigma}=\frac{\partial}{\partial n}v'_{\sigma} &\quad\text{on any Neumann side }\W_1;\\ 
\eta_{\sigma}=v'_{\sigma} &\quad\text{on any Dirichlet side } \W_2;\\ 
\frac{\partial}{\partial n}\eta_{\sigma}=0 &\quad\text{on }S.
\end{dcases}
\end{equation}
Indeed, solutions to such mixed problems are known to exist even on a larger class of Lipschitz domains \cite{brown}. The paper applies to our setting since all angles are strictly less than $\pi$, see the discussion in \cite[Introduction]{brown}. Specifically, since we have at least one Dirichlet side, we may use \cite[Theorem 2.1]{brown} to deduce that 
a solution $\eta_{\sigma}$ to \eqref{eq:mixeddefect} exists, is unique, and 
\[
\|\nabla\eta_{\sigma}\|^2_{L^2(\partial\Omega)}\leq C\left(\left\|\frac{\partial}{\partial n}v'_{\sigma}\right\|^2_{L^2(\W_1)}+\|\nabla_{\tan}v'_{\sigma}\|^2_{L^2(\W_2)}+\|v'_{\sigma}\|^2_{L^2(\W_2)}\right),
\]
where $\nabla_{\tan}$ denotes a tangential derivative along $\W$. Since the Neumann-to-Dirichlet operator on $L^2(\partial\Omega)$ is bounded, again by \cite{Rondi}, we have
\[\|\eta_{\sigma}\|^2_{L^2(\partial\Omega)}\leq C\left\|\frac{\partial}{\partial n}\eta_{\sigma}\right\|^2_{L^2(\partial\Omega)}\leq C\|\nabla\eta_{\sigma}\|^2_{L^2(\partial\Omega)}\]
and therefore
\[
\|\nabla\eta_{\sigma}\|^2_{L^2(\partial\Omega)}+\|\eta_{\sigma}\|^2_{L^2(\partial\Omega)}\leq C\left(\left\|\frac{\partial}{\partial n}v'_{\sigma}\right\|^2_{L^2(\W_1)}+\|\nabla_{\tan}v'_{\sigma}\|^2_{L^2(\W_2)}+\|v'_{\sigma}\|^2_{L^2(\W_2)}\right).
\]

By Theorem \eqref{wedgemodels}, the first two terms on the right hand side are actually bounded by $(C\sigma^{-\mu_\alpha-1})^2$, and the third term is bounded by $(C\sigma^{-\mu_\alpha})^2$. Overall, we obtain
\[
\|\nabla\eta_{\sigma}\|^2_{L^2(\partial\Omega)}+\|\eta_{\sigma}\|^2_{L^2(\partial\Omega)}\leq (C\sigma^{-\mu_\alpha})^2,
\]
so in particular
\begin{equation}\label{eq:etabound}
\|\eta_{\sigma}\|_{L^2(\partial\Omega)}\leq C\sigma^{-\mu_\alpha}.
\end{equation}
The analysis then proceeds as above, and in particular the analogue of Lemma \ref{lem:keyquasi} holds by an identical proof.
\end{remark}

Going back to the Neumann setting again, we let
\begin{equation}\label{quasievals}
\sigma_j=\frac 1L\left(\pi\left(j-\frac 12\right)-\frac{\pi^2}{8}\left(\frac 1\alpha+\frac 1\beta\right)\right),\quad j=1,2,\dots,
\end{equation}
as in \eqref{eq:quantizationcondition}. Abusing notation slightly, let $v_j:=v_{\sigma_j}|_S$ be the traces of the corresponding quasimodes. Assume further that $\|v_j\|_{L^2(S)}=1$. By \eqref{quasievals}. $\sigma_j\geq cj$ for some $c>0$, and thus we have
\begin{equation}\label{keyquasiestrev}
\|\mathcal{D}v_{j}-\sigma_j v_{j}\|_{L^2(S)}\leq Cj^{\delta},\quad \delta=1-\mu_\alpha>0.
\end{equation} 

Now let $\{\varphi_k\}_{k=1}^\infty$ be an orthonormal basis of the eigenfunctions of the sloshing operator $\mathcal{D}$, with eigenvalues $\lambda_j$. By completeness and orthonormality of the $\{\varphi_k\}$, we have, for each $j$,
\begin{equation}\label{eq:fourier}
v_j=\sum_{k=1}^{\infty}a_{jk}\varphi_k,\quad a_{jk}=(v_j,\varphi_k),\quad \sum_{k=1}^{\infty} a_{jk}^2=1.
\end{equation}
Plugging in \eqref{keyquasiestrev}, we get
\[
\left\|\sum_{k=1}^{\infty}a_{jk}(\lambda_k-\sigma_j)\varphi_k\right\|_{L^2(S)}\leq Cj^{-\delta}
\]
and hence
\[
\sum_{k=1}^{\infty}a_{jk}^2(\lambda_k-\sigma_j)^2\leq Cj^{-2\delta}.
\]
Since $\displaystyle\sum_{k=1}^{\infty}a_{jk}^2=1$, it cannot be true that $|\lambda_k-\sigma_j|>Cj^{-\delta}$ for all $k$. Therefore the following Lemma holds.

\begin{lemma}\label{lem:firstlinalgfact} For each $j\in\mathbb N_0$, there exists $k\in\mathbb N_0$ such that
\[
|\sigma_j-\lambda_k|\leq Cj^{-\delta}.
\]
\end{lemma}

In other words, there exists a subsequence of the spectrum that behaves asymptotically as $\sigma_j$, up to an error which is $O(j^{-\delta})$. The key question that we now face is how to prove that $j=k$, i.e., that the sequence gives the {\it  full} spectrum.
In order to achieve this, we have to deal with several issues. First, we do not know whether the quasimodes $\{v_j\}$ form a basis for $L^2(S)$. Second, we do not have good control of the errors $R_{A,\sigma}$ and $R_{B,\sigma}$ near their respective corners. Therefore, some new ideas are needed; in particular, we need more accurate quasimodes. Luckily for us, as will be shown later, we can get away with constructing such quasimodes for some angles only.

\section{Exponentially accurate quasimodes for angles of the form $\frac{\pi}{2q}$}\label{sec:expquasimodes}
\subsection{Hanson-Lewy solutions for angles $\pi/2q$} We will now construct quasimodes in the special case where both angles $\alpha$ and $\beta$ are equal to $\pi/2q$ for some integer $q\geq 2$.

To do this, let us first go back to the study of waves on sloping beaches. Recall that $\Sct{\alpha}$ is the planar sector with $-\alpha\leq\theta\leq 0$, and let $I_1$ be the half-line $\theta=0$ with $I_2$ the half-line $\theta=-\alpha$. Consider the complex variable $z$. Let $q$ be a positive integer and let $\alpha=\pi/2q$, 
\begin{equation}\label{eq:xi}
\xi=\xi_q:=\er^{-\ir \pi/q}.
\end{equation}

\begin{proposition}\label{prop:claim1} Suppose $g(z)$ is an arbitrary function. Set $(\mathcal{A}g)(z):=g(\xi\bar{z})$. Then
\[
(g-\mathcal{A}g)|_{I_2}=0,\quad \frac{\partial}{\partial n}(g+\mathcal{A}g)|_{I_2}=0.
\]
\end{proposition}
\begin{proof} This can be checked by direct computation.
\end{proof}

The following definition is useful throughout.
\begin{definition} The \emph{Steklov defect} of an arbitrary function $g(z)$ on $\Sct{\alpha}$ is
\[
\SD(g):=\left.\left(\frac{\partial g}{\partial y}-g\right)\right|_{I_1}.
\]
\end{definition}
Note that $g(z)$ satisfies the Steklov condition with parameter 1 if and only if $\SD(g)$ is identically zero.

For setup purposes, consider functions of the form
\[
g(z)=\er^{p\bar z},\ h(z)=\er^{pz}.
\]
Differentiation immediately tells us that
\[\
\SD(g)=(-\ir p-1)g|_{I_1};\quad \SD(h)=(\ir p-1)h|_{I_1}.
\]
Given $g=\er^{p\bar z}$, set
\[
\mathcal{B}g(z)=\frac{\ir p+1}{\ir p-1}\er^{pz}.
\]
It is immediate that
\begin{proposition}\label{prop:claim2} For any $g=\er^{p\bar z}$, with $p\in\mathbb C$,
\[
\SD(g+\mathcal{B}g)=0.
\]
\end{proposition}

Now let $f(z)=\er^{-\ir z}$, and consider
\[
f_0(z)=f(z);\ f_1(z)=\mathcal{A}f_0(z),\ f_2=\mathcal{B}f_1(z),\ f_3=\mathcal{A}f_2(z),\dots.
\]
We then have
\[
f_1(z)=f(\xi\bar z),\ f_2(z)=\eta(\xi)f(\xi z),
\]
where
\[
\eta(\xi)=\frac{\ir(-\ir\xi)+1}{\ir(-\ir\xi)-1}=\frac{\xi+1}{\xi-1}.
\]
Further on,
\[
\begin{split}
f_3(z)&=\eta(\xi)f(\xi^2\bar z),\quad f_4(z)=\eta(\xi)\eta(\xi^2)f(\xi^2z),\dots,\\
f_{2q-1}(z)&=\eta(\xi)\eta(\xi^2)\dots\eta(\xi^{q-1})f(\xi^q\bar z)=\prod_{j=1}^{q-1}\eta(\xi^j)f(-\bar z).
\end{split}
\]
Finally, set
\begin{equation}\label{defofualpha}
\upsilon_{\alpha}(z)=\sum_{m=0}^{2q-1}f_m(z).
\end{equation}
\begin{theorem} The function $\upsilon_{\alpha}(z)$ defined by \eqref{defofualpha} is harmonic, satisfies Neumann boundary conditions on $I_2$, and $\SD(\upsilon_{\alpha})=0$.
\end{theorem}
\begin{proof} Since $\upsilon_{\alpha}(z)$ is a sum of rotated plane waves, it is harmonic. It satisfies Neumann boundary conditions because we can write
\[
\begin{split}
\upsilon_{\alpha}&=(f_0+f_1)+(f_2+f_3)+\dots+(f_{2q-2}+f_{2q-1})\\
&=(f_0+\mathcal{A}f_0)+(f_2+\mathcal{A}f_2)+\dots+(f_{2q-2}+\mathcal{A}f_{2q-2})
\end{split}
\]
and use Proposition \ref{prop:claim1} on each term. To see that it satisfies $\SD(\upsilon_{\alpha})=0$, write
\[
\begin{split}
\upsilon_{\alpha}&=f_0+(f_1+f_2)+\dots+(f_{2q-3}+f_{2q-2})+f_{2q-1}\\
&=f_0+(f_1+\mathcal{B}f_1)+\dots+(f_{2q-3}+\mathcal{B}f_{2q-3})+f_{2q-1}.
\end{split}
\]
Using Proposition \ref{prop:claim2} and the linearity of $\SD$, we have
\[
\SD(\upsilon_{\alpha})=\SD(f_0)+\SD(f_{2q-1}),
\]
and all that remains is to show these two terms sum to zero. In fact both are separately zero, because $f_0(z)=\er^{-\ir z}=\er^{-\ir x}\er^y$, and $\SD(f_0)$ is zero by direct computation. Moreover,
\[
f_{2q-1}(z)=\prod_{j=1}^{q-1}\eta(\xi^j)\er^{\ir \bar z}=\prod_{j=1}^{q-1}\eta(\xi^j)\er^{\ir x}\er^y,
\]
and the same direct computation shows $\SD(f_{2q-1})=0$. This completes the proof.
\end{proof}

We call $\upsilon_{\alpha}(z)$ the \emph{Hanson-Lewy solution} for the sloping beach problem with angle $\alpha=\pi/2q$.

It is helpful to introduce the notation
\begin{equation}
\gamma(\xi):=\prod_{j=1}^{q-1}\eta(\xi^j).
\end{equation}
\begin{lemma} We have $\gamma(\xi)=\er^{\ir \pi(q-1)/2}$.
If $q$ is even, $\gamma(\xi)=\pm\ir$, and if $q$ is odd, $\gamma(\xi)=\pm 1$.
\end{lemma}

\begin{proof}We have
\[
\gamma(\xi)=\prod_{j=1}^{q-1}\frac{\er^{-\ir \pi j/q}+1}{\er^{-\ir \pi j/q}-1}=\prod_{j=1}^{q-1}\frac{1+\er^{\ir \pi j/q}}{1-\er^{\ir \pi j/q}},
\]
where we have multiplied numerator and denominator by $\er^{\ir \pi j/q}$. Re-labeling terms in the numerator by switching $j$ for $q-j$ gives
\[
\gamma(\xi)=\prod_{j=1}^{q-1}\frac{1+\er^{\ir \pi (q-j)/q}}{1-\er^{\ir \pi j/q}}=\prod_{j=1}^{q-1}\frac{1-\er^{-\ir \pi j/q}}{1-\er^{\ir \pi j/q}}=\prod_{j=1}^{q-1}\er^{-\ir \pi j/q}\frac{\er^{\ir \pi j/q}-1}{1-\er^{\ir \pi j/q}}=\prod_{j=1}^{q-1}(-\er^{-\ir \pi j/q}).
\]
This may be rewritten as
\[
\gamma(\xi)=\prod_{j=1}^{q-1}(\er^{\ir \pi-\ir\pi j/q})=\exp\left(\ir\pi\left(q-1-\sum_{j=1}^{q-1}\frac jq\right)\right)=\exp\left(\ir\pi\left(q-1-\frac{q-1}{2}\right)\right),
\]
which yields the result.
\end{proof} 

\begin{lemma}\label{lem:decayonrealline} On $I_1$, the Hanson-Lewy solutions $\upsilon_{\alpha}(z)$ are of the form
\[
\upsilon_{\alpha}(x)=\er^{-\ir x}+\gamma(\xi)\er^{\ir x}+\text{ decaying exponentials}.
\]
\end{lemma}
\begin{proof}Indeed, the first and last terms are $\er^{-\ir z}+\gamma(\xi)\er^{\ir \bar z}$, and all other exponentials in the sum defining $\upsilon_{\alpha}$ are of the form $\er^{-\ir \xi^jz}$ or $\er^{-\ir \xi^j\bar z}$ for $j=1,\dots,q-1$. On $I_1$, $z=\bar z=x$, so we have $\er^{-\ir \Re(\xi^j)x}\er^{\Im(\xi^j)x}$. But $\Im(\xi^j)=\sin(-\pi j/q)<0$ for $j=1,\dots,q-1$, and the exponential is thus decaying.
\end{proof}

This may be strengthened: 
\begin{lemma}The rescaled solution $\upsilon_{\alpha}(\sigma z)$ is $\er^{\ir \sigma z}+\gamma(\xi)\er^{\ir \sigma \bar z}$ plus a remainder which is exponentially decaying in $\sigma$ as $\sigma\to\infty$ for each $z\in \Sct{\alpha}$. The exponential decay constant is uniform over all $z\in \Sct{\alpha}$ with $|z|=1$.
\end{lemma}

\begin{proof}  Each term in the sum defining $\upsilon_{\alpha}(z)$ other than the first and last terms is $\er^{-\ir \xi^j z}$ or $\er^{-\ir \xi^j\bar z}$ for $j=1,\dots,q-1$. Observe that $|\er^{-\ir \xi^j\sigma z}|$ decays exponentially in $\sigma$ when $\xi^j z$ is in the negative half-plane, uniformly for $z$ away from the real axis. Since $\arg \xi^j$ is $-\pi j/q$ for $1\leq j\leq q-1$ and $\arg z$ is between $0$ and $-\pi/2q$, we see that for each $z\in \Sct{\alpha}$,
\[
-\frac{\pi}q\leq \arg(\xi^j z)\leq-\pi+\frac{\pi}{2q}.
\]
Thus we have the required decay and it is uniform over the set of $z\in \Sct{\alpha}$ with $|z|=1$. A similar calculation shows that
\[
-\frac{\pi}{2q}\leq \arg(\xi^j\bar z)\leq -\pi+\frac{\pi}{q},
\]
and the same argument applies, completing the proof.
\end{proof}

\begin{remark} A similar construction holds for the solutions of the mixed Steklov-Dirichlet problem.  Set
\[
\tilde f_0=f;\ \tilde f_1=-\mathcal A\tilde f_0;\ \tilde f_2=\mathcal B\tilde f_1;\ \tilde f_3=-\mathcal A\tilde f_2,\dots.
\]
Then set
\[
\tilde \upsilon_{\alpha}(z)=\sum_{e=0}^{2q-1}\tilde f_e(z).
\]
By Propositions \ref{prop:claim1} and \ref{prop:claim2}, these satisfy $\tilde \upsilon_{\alpha}|_{L^2}=0$ and $\SD(u)=0$. They also have the same exponential decay properties.
\end{remark}

Here is a key, novel, observation about the Hanson-Lewy solutions:
\begin{theorem}\label{thm:derivsofhlsols} For each $\alpha=\frac{\pi}{2q}$, the solutions $\upsilon_{\alpha}(z)$ and $\tilde \upsilon_{\alpha}(z)$ satisfy the following properties when restricted to $I_1$:
\[
\upsilon_{\alpha}^{(m)}(0)=0,\ m=q,\dots,2q-1;\quad \tilde \upsilon_{\alpha}^{(m)}(0)=0,\ m=0,\dots,q-1.
\]
\end{theorem}
\begin{proof} Consider first the Neumann case. Set $y=0$. We have
\[
\upsilon_{\alpha}(x)=\er^{-\ir x}+\er^{-\ir \xi x}+\eta(\xi)\er^{-\ir \xi x}+\eta(\xi)\er^{-\ir \xi^2 x}+\dots.
\]
Therefore
\begin{equation}
\label{prodd}
 \upsilon_{\alpha}^{(m)}(0)=(-\ir)^m(1+\xi^m)\left(1+\sum_{k=1}^{q-1}\xi^{km}\prod_{j=1}^k\eta(\xi^j)\right).
\end{equation}
We now apply the Lemma in \cite[p.745]{lewy}. In the notation of Lewy, $\epsilon=\xi=\er^{-\ir \pi/q}$, and the product in the definition of $f(\xi)$ in \cite{lewy} is precisely the one appearing in \eqref{prodd}  times $(-1)^k$. Hence, the right hand side of \eqref{prodd}  vanishes precisely for $m=q,\dots,2q-1$. Indeed, for $m=q$ this is obvious since $1+\xi^q=0$; for $m=q+r$, $r=1,\dots,q-1$, we have $\xi^{k(q+r)}=\xi^r(\xi^q)^k=(-1)^k\xi^{k r}$, and the result follows from Lewy's lemma. 

In the Dirichlet case we have
\[
\tilde \upsilon_{\alpha}(x)=\er^{-\ir x}-\er^{-\ir \xi x}-\eta(\xi)\er^{-\ir \xi x}+\eta(\xi)\er^{-\ir \xi^2 x}+\dots.
\]
Therefore
\[
\tilde \upsilon_{\alpha}^{(m)}(0)=(-\ir^m)(1-\xi^m)\left(1+\sum_{k=1}^{q-1}(-1)^k\xi^{km}\prod_{j=1}^{k}\eta(\xi^j)\right).
\]
A similar use of Lewy's lemma shows this is zero for $m=0,1,\dots,q-1$, and the theorem is proved. 
\end{proof}
\begin{remark} An alternative proof of this Theorem was provided to us by Rinat Kashaev \cite{Kashaev}, based on combinatorial techniques used in \cite{Kashaev1}.
\end{remark}

\begin{remark}
 The Hanson-Lewy solutions are closely connected to a higher order Sturm-Liouville type ODE on the real line, specifically
\begin{equation}\label{eq:plainODE}
(-1)^qU^{(2q)}=\Lambda^{2q}U.
\end{equation}
It is immediate that the Hanson-Lewy solutions $\upsilon_{\alpha}$ and $\tilde \upsilon_{\alpha}$ satisfy this equation with $\Lambda=1$. This observation together with Theorem \ref{thm:derivsofhlsols} will be key for the next step of the argument.
\end{remark}

\subsection{Exponentially accurate quasimodes for angles of the form $\pi/2q$} 
\label{expacc}
Consider the sloshing problem on a triangle with top side of length $L=1$ (for simplicity) and with angles $\alpha=\beta=\pi/2q$, $q\in\mathbb N$, $q$ even. Suppose that $\sigma$ satisfies the quantization condition \eqref{eq:quantizationcondition}. For such $\sigma$, let
\[
g_{\sigma}(z)=\upsilon_{\alpha}(\sigma z)+\upsilon_{\beta}^\text{d}(\sigma(1-\overline{z})),
\]
where $\upsilon_{\beta}^\text{d}$ is the sum of just the \emph{decaying} exponents in $\upsilon_{\beta}$ (see Lemma \ref{lem:decayonrealline}). Note that the choice of $\sigma$ yields that the oscillating exponents in $\left(\upsilon_{\alpha}|_{I_1}\right)(\sigma x)$ and $\left(\upsilon_{\beta}|_{I_1}\right)(\sigma(1-x))$ coincide, so there is a similar decomposition near $B$.

First we view the functions $g_{\sigma}$ as quasimodes for the sloshing operator. We are in the same situation we were with the Peters solutions, but now the errors decay exponentially rather than  polynomially. By an identical argument to that in subsection \ref{subsection:step2}, there exist functions $\eta_{\sigma}(z)$ satisfying \eqref{eq:etabound} such that if we set $v_{\sigma}(z)=g_{\sigma}(z)-\eta_{\sigma}(z)$, then we have the key quasimode estimate
\begin{equation}\label{eq:gquasiforsloshing}
\|\mathcal{D}(v_{\sigma}|_S)-\sigma v_{\sigma}\|_{L^2(S)}=O\left(\er^{- \sigma/C}\right).
\end{equation}

Now we view $g_{\sigma}$ as quasimodes for an ODE problem. By Theorem \ref{thm:derivsofhlsols}, we see that $g_{\sigma}(x)$ is a solution of the ODE \eqref{eq:plainODE}, and moreover one which satisfies up to an error $O(\er^{-\sigma/C})$, for some $C>0$, the self-adjoint boundary conditions
\[
g_{\sigma}^{(m)}(0)=g_{\sigma}^{(m)}(1)=0,\quad m=q,\dots,2q-1.
\]
at the ends of the interval. Specifically, we claim that these functions $g_{\sigma}(x)$ are exponentially accurate quasimodes  on $[0,1]$ for the elliptic self-adjoint ODE eigenvalue problem
\begin{equation}\label{odeneumann}
\begin{dcases}
(-1)^qU^{(2q)}=\Lambda^{2q}U\\ 
U^{(m)}(0)=U^{(m)}(1)=0\text{ for }m=q,\dots,2q-1.
\end{dcases}
\end{equation}
In order to show this, we need to run a similar argument to the one we have run in subsection \ref{subsection:step2} to correct our quasimodes. First, we correct the boundary values of $g_{\sigma}$ in order for the quasimode to lie in the domain of the operator. We do it as follows.

Let $\Phi_m(x)$ be a function equal to $x^m$ near $x=0$ and smoothly decaying so that it is identically zero whenever $x>1/2$. Then for each $\sigma=\sigma_k$, there exists a function 
\[
\bar\eta_{\sigma}(x)=\sum_{m=q}^{2q-1}\left(a_m\Phi_m(x)+b_m\Phi_m(1-x)\right),
\]
where and $a_m$, $b_m$ are both $O\left(\er^{-\sigma/C}\right)$, and where $g_{\sigma}(x)+\bar\eta_{\sigma}(x)$ satisfy the boundary conditions in \eqref{odeneumann}. In other words,
\[
\bar v_{\sigma}(x)=g_{\sigma}(x)+\bar\eta_{\sigma}(x)
\]
could be used as a quasimode for the eigenvalue problem \eqref{odeneumann}. Indeed, it is easy to see that
\begin{equation}\label{eq:gquasiforode}
\|(-1)^q\bar v_{\sigma}^{(2q)}(x)-\sigma^{2q}\bar v_{\sigma}(x)\|_{L^2([0,1])}=O\left(\er^{-\sigma/C}\right).
\end{equation}
The consequences of these quasimode estimates will be discussed in the next section. Note that by our estimates on $\eta$, the functions $\tilde\eta$, $v_{\sigma}$, $g_{\sigma}$, and $\tilde v_{\sigma}$ are all within $O\left(\er^{-\sigma/C}\right)$ in $L^2(S)$ norm, and in particular 
\begin{equation}\label{eq:moralequivalence}
||\bar v_{\sigma}-v_{\sigma}||_{L^2(S)}\leq C\er^{-\sigma/C}.
\end{equation}

\section{Completeness of quasimodes via higher order Sturm-Liouville equations}\label{sec:completeness}
\subsection{An abstract linear algebra result} In order to complete the proof of Theorem \ref{thm:main} for triangular domains, we need to show that the  family of quasimodes constructed in subsection \ref{subsection:step2} is complete, and make sure that the numeration of the corresponding eigenvalues is precisely as in \eqref{eq:sloshingasymptotics}.
 Unlike the general case when we were using Peters solutions, we are well equipped to do that for angles $\pi/2q$ using exponentially accurate quasimodes. We will need the following strengthening of Lemma \ref{lem:firstlinalgfact}, formulated as an abstract linear algebra result below.

\begin{theorem}\label{thm:linalg} Let $\boldsymbol{D}$ be a self-adjoint operator on an infinite-dimensional Hilbert space $\boldsymbol{H}$ with a discrete spectrum $\{\lambda_j\}_{j=1}^\infty$ and a complete orthonormal basis of corresponding eigenvectors $\{\boldsymbol{\varphi}_j\}$. Suppose $\{\boldsymbol{v}_j\}$ is a sequence of  quasimodes (with  $\|\boldsymbol{v}_j\|_{\boldsymbol H}=1$) such that
\begin{equation}\label{thatestimate}
\|\boldsymbol D \boldsymbol v_j-\sigma_j \boldsymbol v_j\|_{\boldsymbol H}\leq F(j),
\end{equation}
where  $\lim_{j\to\infty} F(j)=0$. Then
\begin{enumerate}
\item[1)] For all $j$, there exists $k$ such that $|\sigma_j-\lambda_k|\leq F(j)$;
\item[2)] For each $j$, there exists a vector $\boldsymbol{w}_j$ with the following properties:
\begin{itemize}
\item $\boldsymbol{w}_j$ is a linear combination of eigenvectors of $\mathcal D$ with eigenvalues in the interval $[\sigma_j-\sqrt{F(j)},\sigma_j+\sqrt{F(j)}]$;
\item $\|\boldsymbol{w}_j\|_{\boldsymbol{H}}=1$;
\item $\|\boldsymbol{w}_j-\boldsymbol{v}_j\|_{\boldsymbol{H}}\leq\sqrt{F(j)}+\sqrt{\frac{F(j)}{1-F(j)}}= 2\sqrt{F(j)}(1+o(1))$ as $j\to\infty$.
\end{itemize}
\end{enumerate}
\end{theorem}
\begin{proof} Indeed, the first part of the statement is proved  using precisely the same arguments as  Lemma \ref{lem:firstlinalgfact}. Let us prove part 2) of the theorem.  Set, similarly to \eqref{eq:fourier}, 
\[
\boldsymbol v_j=\sum_{k=1}^{\infty}\boldsymbol a_{jk}\boldsymbol \varphi_k,\quad \boldsymbol a_{jk}=(\boldsymbol v_j,\boldsymbol \varphi_k)_{\boldsymbol H},\quad \sum_{k=0}^{\infty} \boldsymbol a_{jk}^2=1.
\]
From the estimate \eqref{thatestimate},
\begin{equation}
F(j)^2\geq \sum_{k:\ |\lambda_k-\sigma_j|>\sqrt{F(j)}} \boldsymbol a_{jk}^2(\lambda_k-\sigma_j)^2\geq F(j)\sum_{k: |\lambda_k-\sigma_j|>\sqrt{F(j)}} \boldsymbol a_{jk}^2,
\end{equation}
hence
\begin{equation}
\sum\limits_{k:\ |\lambda_k-\sigma_j|>\sqrt{F(j)}} \boldsymbol a_{jk}^2\leq F(j)\Rightarrow 1-F(j)\leq \sum\limits_{k:\ |\lambda_k-\sigma_j|\leq \sqrt{F(j)}} \boldsymbol a_{jk}^2\leq 1.
\end{equation}
Letting
\[
\boldsymbol w_j:=\dfrac{\sum\limits_{k:\ |\lambda_k-\sigma_j|\leq \sqrt{F(j)}}\boldsymbol a_{jk}\boldsymbol \varphi_k}{\sum\limits_{k:\ |\lambda_k-\sigma_j|\leq \sqrt{F(j)}} \boldsymbol a_{jk}^2}
\]
completes the proof.
\end{proof}
As a consequence of \eqref{eq:gquasiforsloshing} and \eqref{eq:gquasiforode}, this Lemma may be applied to both the sloshing and ODE problems in the case where both $\alpha$ and $\beta$ equal $\pi/2q$. In that case, for each $j$ sufficiently large, there exist functions $w_j$ and $W_j$ which are both within $O\left(\er^{-j/C}\right)$ of $v_j$ in $L^2$ norm (by \eqref{eq:moralequivalence}), and which are linear combinations of eigenfunctions for the sloshing and ODE problems, respectively, with eigenvalues within $O\left(\er^{-j/C}\right)$ of $\sigma_j$. Therefore there is an infinite subsequence, also denoted $w_j$, of linear combinations of sloshing eigenfunctions  which are close to linear combinations of eigenfunctions $W_j$ of the ODE. The problem now is to prove asymptotic completeness and to make sure that the enumeration lines up.

\subsection{Eigenvalue asymptotics for  higher order Sturm-Liouville problems}\label{naimark}
Our next goal is to understand the asymptotics of the eigenvalues of the higher order Sturm-Liouville problem \eqref{odeneumann} with Neumann boundary conditions.
For  Dirichlet boundary conditions, the corresponding eigenvalue problem is 
\begin{equation}\label{odedirichlet}
\begin{dcases}(-1)^q U^{(2q)}=\Lambda^{2q}U\\ U^{(m)}(0)=U^{(m)}(1)=0, m=0,\dots,q-1.\end{dcases}
\end{equation}
As was mentioned in subsection \ref{stliouv}, the Dirichlet and Neumann spectra for the ODE above are related by the following proposition: 
\begin{proposition}\label{prop:weirdODEfact} The nonzero eigenvalues of the Neumann ODE problem \eqref{odeneumann} are the same as those of the Dirichlet ODE problem \eqref{odedirichlet}, including multiplicity. However, the kernel of \eqref{odedirichlet} is trivial, whereas the kernel of \eqref{odeneumann} consists of polynomials of degree less than $q$ and therefore has dimension $q$.
\end{proposition}
\begin{proof} The kernel claims follow from direct computation. For the rest, the solutions with nonzero eigenvalue of the ODE problems \eqref{odeneumann} and \eqref{odedirichlet} must be of the form
\begin{equation}\label{eq:genodesol}
U(x)=\sum_{k=0}^{2q-1}c_k\er^{\omega_k\Lambda x},
\end{equation}
where $\{\omega_1,\dots,\omega_n\}$ are the $n$th roots of $-1$. We claim that a particular $U(x)$ satisfies the Dirichlet boundary conditions if and only if the function
\[
\tilde U(x):=\sum_{k=0}^{2q-1}c_k\omega_k^q\er^{\omega_k\Lambda x}
\]
satisfies the Neumann boundary conditions. Indeed, this is obvious; since $\omega_k^{2q}=1$, for any $m$ and any $x$, $U^{(m)}(x)=0$ if and only if $\tilde U^{(m+q)}(x)=0$. This gives a one-to-one correspondence between Dirichlet and Neumann eigenfunctions which completes the proof.
\end{proof}
\begin{remark}
Proposition \ref{prop:weirdODEfact} also follows from  the observation that the operators corresponding to the problems
\eqref{odeneumann} and \eqref{odedirichlet}  can be represented as $AA^*$ and $A^*A$, where $A$ is an operator given by 
$A\, U=\ir^q U^{(q)}$ subject to boundary conditions $U^{(m)}(0)=U^{(m)}(1)=0, m=0,\dots,q-1.$
\end{remark}
The following asymptotic result may be extracted, with some extra work (see Appendix \ref{naimarkproof}), from a  book of M.~Naimark  \cite[Theorem 2, equations (45 a) and (45 b)]{naimark}. Note that our differential operator has order $n=2q$, so $\mu$ in the notation of \cite{naimark} is equal to $q$. 
\begin{proposition}\label{prop:odeasymp} For each $q\in \mathbb{N}$, all sufficiently large eigenvalues of the boundary value problem  \eqref{odeneumann} are given by the formulae 
\begin{equation}\label{odeasymptotics}
\begin{split}
\left\{\left(k-\frac 12\right)\pi+O\left(\frac 1k\right)\right\}_{k=K}^\infty\qquad&\text{for $q$ even},\\
\\
\left\{k\pi+O\left(\frac 1k\right)\right\}_{k=K}^\infty\qquad&\text{for $q$ odd},
\end{split}
\end{equation}
for some $K\in\mathbb{N}$.
\end{proposition}

\begin{corollary}\label{cor:odeasymp}
In particular, there exists a $J=J_q\in\mathbb Z$ so that for large enough $j$,
\[
\Lambda_{j}=\sigma_{j+J}+ O\left(\frac 1j\right).
\]
\end{corollary}
This corollary is immediate from the explicit formula for $\sigma_j$, with the appropriate values of $\alpha=\beta=\pi/2q$.

\begin{remark} Proposition \ref{prop:odeasymp}  implies that $\lambda_j$ are simple and separated by nearly $\pi$ for large enough $j$. However, it does not  tell us right away that $J=0$,
and further work is needed to establish this.  \end{remark}

The proof of Proposition \ref{prop:odeasymp} for even $q$ is given in  Appendix \ref{naimarkproof}. 
\begin{remark}
\label{naimark:proof}
The argument presented in \cite{naimark} is slightly different for $q$ even and $q$ odd. In Appendix \ref{naimarkproof} we shall assume that $q$ is even, which is in fact sufficient for our purposes. The proof for odd $q$ is analogous and is left to the interested reader. 
\end{remark}

\subsection{Connection to sloshing eigenvalues} Recall now that $\{w_j\}$ are linear combinations of sloshing eigenfunctions and $\{W_j\}$ are linear combinations of ODE eigenfunctions, each with eigenvalues in shrinking intervals around $\sigma_j$, and each exponentially close to a quasimode. As a consequence of the remark following Corollary \ref{cor:odeasymp}, there are gaps between consecutive eigenvalues of the ODE and so for large enough $j$ there is only one eigenvalue of the ODE in each of those intervals. This means that $W_j$, for sufficiently large $j$, must actually be an eigenfunction of the ODE, rather than a linear combination of eigenfunctions. By \ref{cor:odeasymp}, we must have $W_j=U_{j-J}$, with eigenvalue $\Lambda_{j-J}$.

Now consider $\{w_j\}_{j\geq N}$ and $\{W_j\}_{j\geq N}$, and let their spans be $\boldsymbol X_N$ and $\boldsymbol X_N^*$ respectively. Pick $N$ large enough so that $W_j$ are eigenfunctions and so that the intervals in Theorem \ref{thm:linalg} are disjoint, and also large enough so that
\begin{equation}\label{eq:bigenoughN}
\sum_{j\geq N}\|w_j-W_j\|^2_{L^2(S)}\leq\sum_{j\geq N}\left(C\er^{-j/C}\right)^2<1.
\end{equation}
We know by ODE theory that $\{W_j\}$ form a complete orthonormal basis of $L^2(S)$ and that therefore \[\dim(\boldsymbol X_N^*)^{\perp}=N-1.\]
\begin{lemma}\label{lem:barykrein} For sufficiently large $N$, the dimension $\dim \boldsymbol X_N^{\perp}=N-1$ as well.
\end{lemma}
\begin{proof} This is essentially a version of the Bary-Krein lemma (see \cite{Kazdan}) and our argument closely follows \cite{stackexchange}. Define $\boldsymbol{A}:\boldsymbol X_N^{\perp}\to(\boldsymbol X_N^*)^{\perp}$ and $\boldsymbol{B}:(\boldsymbol X_N^*)^{\perp}\to \boldsymbol X_N^{\perp}$ by
\[
\boldsymbol{A}f=f-\boldsymbol{P}_{\boldsymbol X_N^{*}}f;\ \boldsymbol{B}f=f-\boldsymbol{P}_{\boldsymbol X_N}f,
\]
where $\boldsymbol{P}_{\boldsymbol X_N^*}$ and $\boldsymbol{P}_{\boldsymbol X_N}$ are orthogonal projections. By \eqref{eq:bigenoughN} we have for any $f\in \boldsymbol X_N^{\perp}$:
\[
\|\boldsymbol{P}_{\boldsymbol X_N^*}f\|^2=\sum_{j\geq N}|(f,w_j^*)|^2=\sum_{j\geq N}|(f,w_j^*-w_j)|^2\leq \epsilon\|f\|^2,\ \epsilon<1.
\]
Similarly, for any $f\in (\boldsymbol X_N^*)^{\perp}$, $\|\boldsymbol{P}_{\boldsymbol X_N}f\|<\epsilon^2\|f\|^2$. Therefore, for any $f\in \boldsymbol X_N^{\perp}$,
\[
\boldsymbol{B}\boldsymbol{A}f=(f-\boldsymbol{P}_{\boldsymbol X_N^*}f)-\boldsymbol{P}_{\boldsymbol X_N}(f-\boldsymbol{P}_{\boldsymbol X_N^*}f)=f-(I-\boldsymbol{P}_{\boldsymbol X_N})\boldsymbol{P}_{\boldsymbol X_N^*}f=f-\boldsymbol{P}_{\boldsymbol X_N^{\perp}}\boldsymbol{P}_{\boldsymbol X_N^*}f,
\]
and hence
\[
\|(I_{\boldsymbol X_N^{\perp}}-\boldsymbol{B}\boldsymbol{A})f\|=\|\boldsymbol{P}_{\boldsymbol X_N^{\perp}}\boldsymbol{P}_{\boldsymbol X_N}^*f\|\leq \|\boldsymbol{P}_{\boldsymbol X_N^*}f\|<\epsilon\|f\|\mbox{ for }\epsilon<1.
\]
Therefore, $\boldsymbol{A}$ must be injective, as otherwise $\boldsymbol{B}\boldsymbol{A}f=0$ for some nonzero $f$ and we get a contradiction. Hence $\dim \boldsymbol X_N^{\perp}\leq\dim(\boldsymbol X_N^{*})^{\perp}=N-1$. Repeating the same argument with $\boldsymbol X_N$ and $\boldsymbol X_N^{*}$ interchanged shows that $\dim(\boldsymbol X_N^*)^{\perp}\leq\dim \boldsymbol X_N^{\perp}$, and therefore that
\[
\dim \boldsymbol X_N^{\perp}=\dim (\boldsymbol X_N^*)^{\perp}=N-1,
\]
as desired. \end{proof}

One immediate consequence of this Lemma is that the sequence $\{w_j\}_{j\geq N}$ after a certain point contains only pure eigenfunctions. Indeed, if $w_j$ is a linear combination of $k\geq 2$ eigenfunctions, there are $k-1$ linearly independent functions generated by the same eigenfunctions, orthogonal to $w_j$ (and all other $\{w_j\}$, since eigenfunctions are orthonormal). 
But since $\dim \boldsymbol X_N^{\perp}<\infty$, $k=1$ starting from some $j=N$. Without loss of generality we can pick $N$ sufficiently large so that all $w_j$ are simple for $j\geq N$.

Even more importantly, Lemma \ref{lem:barykrein} tells us that the sequences $\{w_j\}_{j\geq N}$ and $\{W_j\}_{j\geq N}$ are missing the same number of eigenfunctions, namely $N-1$. Since we know that $W_j$ corresponds to eigenvalue $\Lambda_{j-J}=\sigma_j+O\left(\frac{1}{j}\right)$, we have proved the following

\begin{proposition}\label{thm:asympwithshift} In the case $\alpha=\beta=\pi/2q$,  there exists a constant $C>0$ and an integer $J_{q}$ such that both the sloshing eigenvalues 
$\lambda_j$ and the ODE eigenvalues $\Lambda_j$ satisfy
\begin{equation}
\label{withshifts}
\lambda_j L=\pi\left(j+J_q-\frac 12-\frac{q}{2}\right)+ O\left(\er^{-Cj}\right) = \Lambda_j L + O(\er^{-Cj}) ,
\end{equation}
\end{proposition}
\begin{remark}
Under the same assumptions, asymptotics \eqref{withshifts} holds for the Steklov-Dirichlet eigenvalues $\lambda_j^D$  and the eigenvalues $\Lambda_j^D$ of the ODE with Dirichlet boundary conditons \eqref{odedirichlet}  with  $-\frac{q}{2}$ being replaced by  $+\frac{q}{2}$. Note that the shift $J_q$ in the Dirichlet case (which a priori could be different) is the same as in the Neumann case, as immediately follows from Proposition \ref{prop:weirdODEfact}.
\end{remark}

\subsection{Proof of Theorem \ref{ODE} for triangular domains}
We are now in a position to complete the proof of Theorem \ref{ODE}.
Given Proposition \ref{thm:asympwithshift}, it remains  to show that the shift $J_q=0$ for all $q\in\mathbb N$. 
We prove this by induction.  For $q=1$, by computing the eigenvalues of the standard (second-order) Sturm-Liouville problem, we observe that $J_1=0$. 
In order to make the induction step we  use the domain monotonicity properties of Steklov-Neumann and Steklov-Dirichlet eigenfunctions, namely:
\begin{proposition}\cite [section 3]{BKPS}
\label{prop:domainmonotonicity} Suppose $\Omega_1$ and $\Omega_2$ are two domains with the same sloshing surface $S$, with $\Omega_1\subseteq\Omega_2$. Then for all $k=1,2,\dots$, we have $\lambda_k(\Omega_1)\leq\lambda_k(\Omega_2)$ and 
$\lambda_k^D(\Omega_1)\geq\lambda_k^D(\Omega_2)$.
\end{proposition}
Let $\Omega_q$ be an isosceles triangle with angles $\pi/2q$ at the base. We have $\Omega_{q+1} \subset \Omega_q$ for all $k\ge 1$.
Assume now $J_{q+1}>J_q=0$ for some $q>0$. Then it immediately follows from \eqref{withshifts} that $\lambda_j(\Omega_{q+1})>\lambda_j(\Omega_q)$ for large $j$, which contradicts
Proposition \ref{prop:domainmonotonicity} for Neumann eigenvalues. Assuming instead that $J_{q+1}<J_q$, we also get a contradiction, this time with Proposition \ref{prop:domainmonotonicity} for Dirichlet eigenvalues. Therefore, $J_q=0$ for all $k=1,2,\dots$, and this completes the proof of Theorem \ref{ODE} for triangular domains.. \qed
 \begin{remark}  A surprising feature of this proof is that domain monotonicity for mixed Steklov eigenvalues implies new results for the eigenvalues of 
higher order Sturm-Liouville problems, which a priori are easier to investigate.
\end{remark}

\subsection{Completeness of quasimodes for arbitrary triangular domains} 
\label{subsection:monotonicity}

We now complete the proofs  of Theorems  \ref{thm:main} and \ref{thm:Dirichlet}, as well as Propositions \ref{prop:pi2} and \ref{prop:pi22},  for triangular domains, by taking full advantage of the domain monotonicity properties in Proposition \ref{prop:domainmonotonicity}.

Let $\Omega_s$, $s\in[0,1]$, be a continuous family of sloshing domains (not necessarily triangles) sharing a common sloshing surface $S$ with $L=1$. Assume that each $\Omega_s$ is straight in a neighbourhood of the vertices, with angles $\alpha(s)$ and $\beta(s)$. Moreover, assume that $\Omega_s$ is a monotone family, i.e. that $s<t\Rightarrow \Omega_s\subseteq\Omega_t$, and assume that $\alpha(s)$ and $\beta(s)$ are both less than $\pi/2$ for all $0\leq s<1$ (possibly equaling $\pi/2$ for $s=1$).

For the moment we specialise to the Neumann setting; denote the associated quasi-frequencies by $\sigma_j(s)$, and observe by the formula \eqref{eq:quantizationcondition} that they are uniformly equicontinuous in $s$. By Lemma \ref{lem:firstlinalgfact}, for any $s$ and any $j$, we know that there exist integers $k(j,s)$ such that for $j$ sufficiently large,
\begin{equation}\label{eq:windows}
|\sigma_j(s)-\lambda_{k(j,s)}(s)|=o(1). 
\end{equation}
 By the work of Davis \cite{dav65}, we have a similar bound  if $\alpha(s)=\beta(s)=\pi/2$.

The completeness property we need to show  translates to the statement that $k(j,s)=j$ for all sufficiently large $j$, for then the indices match and we have decaying bounds on $|\sigma_j-\lambda_j|$, rather than just $|\sigma_j-\lambda_k|$ for some unknown $k$.

The key lemma is the following.
\begin{lemma}\label{lem:monotonicityconsequences} Suppose that $\Omega_s$, $s\in [0,1]$, is a family of sloshing domains as above. Consider the Neumann case. Then
\begin{enumerate}
\item If $s<s'$ and $k(j,s')\geq j$ for all sufficiently large $j$, then there exists $N>0$ so that for all $j\geq N$, $k(j,s)\geq j$.
\item If $s<s'$ and $k(j,s)\leq j$ for all sufficiently large $j$, then there exists $N>0$ so that for all $j\geq N$, $k(j,s')\leq j$.
\item If both $k(j,0)=j$ and $k(j,1)=j$ for all sufficiently large $j$, then for each $s\in [0,1]$, there exists $N>0$ so that for all $j\geq N$, $k(j,s)=j$.
\end{enumerate}
In the Dirichlet case, the same result holds if we flip the inequalities in the conclusion of the first two statements.
\end{lemma}
\begin{remark} Of course, $k(j,s)$ may not be well-defined for small $j$, as there may in some cases be more than one value that works, but if the hypotheses are satisfied for some choices of $k$, then the conclusions must be satisfied for all choices of $k$.
\end{remark}

\begin{proof}[Proof of Lemma \ref{lem:monotonicityconsequences}] The third statement is an immediate consequence of the first two (applied with $s'=1$ and $s=0$ respectively), so we need only to prove the first two. Since the second one is practically identical to the first one, we only prove the first one here.

First consider the case when $s$ and $s'$ are such that $|\sigma_j(s)-\sigma_j(s')|$ (which is independent of $j$) is less than $\pi$, specifically less than $\pi-\epsilon$ for some $\epsilon>0$. Then for sufficiently large $j$, applying Lemma \ref{lem:firstlinalgfact} as in the discussion above, we can arrange both 
\[|\lambda_{k(j,s)}(s)-\lambda_{k(j,s')}(s')|<\pi-\epsilon/2\mbox{ and }\lambda_{k(j-1,s')}(s')<\lambda_{k(j,s')}(s')-(\pi-\epsilon/2).\]
As a consequence,
\[\lambda_{k(j-1,s')}(s')<\lambda_{k(j,s)}(s).\]
But by domain monotonicity applied to $\lambda_{k(j-1,s')}$, this implies that
\[\lambda_{k(j-1,s')}(s)<\lambda_{k(j,s)}(s),\]
and therefore that $k(j,s)>k(j-1,s')$. Since $k(j-1,s')\geq j-1$ for sufficiently large $j$ by assumption, we must have $k(j,s)>j-1$ and hence $k(j,s)\geq j$, since it is an integer. This is what we wanted.

In the case where $s$ and $s'$ are not such that $|\sigma_j(s)-\sigma_j(s')|<\pi$, simply do the proof in steps: first extend to some $s_1<s$ with $|\sigma_j(s_1)-\sigma_j(s')|<\pi$, then to some $s_2<s_1$, et cetera. Since in all cases $|\sigma_j(s)-\sigma_j(s')|$ is finite, this process can be set up to terminate in finitely many steps, completing the proof.
\end{proof}

Now we can complete the proofs  of Theorems \ref{thm:main} and \ref{thm:Dirichlet}  for triangular domains. 
It is possible to embed the triangle $\Omega$ as $\Omega_{1/2}$ in a continuous, nested family of sloshing domains $\{\Omega_{s}\}_{s\in[0,1]}$, where $\Omega_0$ is an isosceles triangle with angles $\pi/2q$ for some large even $q$, and where $\Omega_1$ is a \emph{rectangle}. For $\Omega_1$, explicit calculations (see \cite{BKPS}) show that
\[
\lambda_kL=\pi\left(k-\frac 12\right)-\frac{\pi^2}{8}\left(\frac{2}{\pi}+\frac{2}{\pi}\right)+o(1)=\pi(k-1) +o(1),
\]
and therefore that $k(j,1)=j$ for all sufficiently large $j$, which is our completeness property. However, we have already proved completeness for $\Omega_0$, in Theorem \ref{ODE}, so $k(j,0)=j$ for all sufficiently large $j$. And we know we can construct quasimodes on $\Omega_{1/2}$, since it is a triangle, and therefore by Lemma \ref{lem:firstlinalgfact} there exist integers $k(j,1/2)$ such that
\[
\left|\sigma_j-\lambda_{k(j,1/2)}\right|=O\left(k^{1-\frac{\pi/2}{\max\{\alpha,\beta\}}}\right).
\]
By statement (3) of Lemma \ref{lem:monotonicityconsequences}, we have $k(j,1/2)=j$ for all sufficiently large $j$. Thus for all sufficiently large $k$,
\[
\lambda_kL=\sigma_kL+O\left(k^{1-\frac{\pi/2}{\max\{\alpha,\beta\}}}\right)=\pi\left(k-\frac 12\right)-\frac{\pi^2}{8}\left(\frac 1{\alpha}+\frac 1{\beta}\right)+O\left(k^{1-\frac{\pi/2}{\max\{\alpha,\beta\}}}\right),
\]
which proves Theorem \ref{thm:main}. The proof of Theorem \ref{thm:Dirichlet} is essentially identical, with domain monotonicity going in the opposite direction, and the estimate on the error term being modified according to Theorem \ref{wedgemodels}. 

The proofs of Propositions \ref{prop:pi2} and \ref{prop:pi22} for triangular domains also follow in the same way, by using the model solutions \eqref{pi2solutions} in the $\Sct{\pi/2}$ in order to construct the quasimodes.

\subsection{Completeness of quasi-frequencies: abstract setting}
\label{compabst}
Here, we would like to formulate the results of the previous subsection on the abstract level. We will not need these results in the present paper, but we will need them
in the subsequent paper \cite{LPPS}.
The proofs of these results are very similar to the proofs above, and we will omit them. 

Let $\Sigma:=\{\sigma_j\}$ and  $\Lambda:=\{\lambda_j\}$, $j=1,2,...$ be two non-decreasing sequences of real numbers tending to infinity, which are called {\it quasi-frequencies} and {\it eigenvalues},  respectively. Suppose that the following ``quasi-frequency gap'' condition is satisfied: there exists a constant $C>0$ such that 
\begin{equation}
\label{qgap}
\sigma_{j+1}-\sigma_j>C
\end{equation}
for sufficiently large $j$.

\begin{definition}
We say that $\Sigma$ is a{\it  system of quasi-frequencies approximating eigenvalues $\Lambda$} 
if there exists a mapping $k:{\mathbb N}\to{\mathbb N}$ 
such that: 

(i)  $|\sigma_j-\lambda_{k(j)}| \to 0$;  

\smallskip

(ii)  $k(j_1)\ne k(j_2)$   for sufficiently large and distinct  $j_1$, $j_2$.
 
\end{definition}

Under assumptions (i-ii) and \eqref{qgap}  we obviously have that for $j$ large enough, $\{k(j)\}$ is a strictly increasing sequence of natural numbers. Therefore,  if we denote $m(j):=k(j)-j$, then the sequence $\{m(j)\}$ is non-decreasing for $j$ large enough,  and, therefore, it is converging (possibly, to $+\infty$). Put
$$
M=M(\Sigma,\Lambda):=\lim_{j\to\infty} m(j). 
$$

\begin{definition}
We say that the system of quasi-frequencies $\Sigma$ is  \emph{asymptotically complete in} $\Lambda$ if $M\ne +\infty$. 
\end{definition}

Now we assume that we have one-parameter families $\Sigma(s)=\{\sigma_{j}(s)\}$ and $\Lambda(s)=\{\lambda_{j}(s)\}$, $s\in [0,1]$ of quasi-frequencies and eigenvalues. 

\begin{lemma}
Suppose the family $\Sigma(s)$ is equicontinuous and the family $\Lambda(s)$ is monotone increasing (i.e. each $\lambda_{j}(s)$ is increasing). Then for sufficiently large $j$ (uniformly in $s$) the function $m(j)=m(j,s)$ is non-increasing in $s$.  
\end{lemma}
\begin{proof}
The proof is the same as the proof of statement (i) of Lemma \ref{lem:monotonicityconsequences}.
\end{proof}
Under the conditions of the previous Lemma, the function  $M=M(s)$ is non-increasing. Therefore, the following alternative formulation of statement (iii) of Lemma \ref{lem:monotonicityconsequences} can be obtained as an immediate consequence:
\begin{corollary}
If \,  $\Sigma(0)$ is asymptotically complete in $\Lambda(0)$, then $\Sigma(1)$ is asymptotically complete in $\Lambda(1)$. 
\end{corollary}
\begin{remark}  With some extra work the ``quasi-frequency'' gap'' assumption \eqref{qgap} could be weakened,  and it is sufficient to require that the number of quasi-frequencies $\sigma_j$ inside any interval of length one is bounded. More details on that will be provided in \cite{LPPS}.

\end{remark}

\subsection{Proof of Proposition \ref{thm:mixed}}  Using a similar strategy, let us now  prove the asymptotics \eqref{mixedasymp} for the mixed Steklov-Dirichlet-Neumann problem on triangles. 

First, assume that $\alpha=\pi/2$. Consider a sloshing problem on the doubled isosceles triangle $\bigtriangleup B'BZ$, where $B'$ is symmetric to $B$ with respect to $A$. Using the odd-even decomposition of eigenfunctions with respect to symmetry across $AZ$, we see that the spectrum of the sloshing problem on $\bigtriangleup B'BZ$ is the union (counting multiplicity) of the eigenvalues of the sloshing problem and of the mixed Steklov-Dirichlet-Neumann problem on $\bigtriangleup ABZ$. Given that we have already computed the asymptotics of the two former problems, it is easy to check that the spectral asymptotics of the latter problem satisfies \eqref{mixedasymp}.  A similar reflection argument works in the case when $\beta=\pi/2$ and $\alpha$ is arbitrary. This proves the proposition when either of the angles is equal to $\pi/2$. 

Let now $0<\alpha,\beta\le \pi/2$ be arbitrary. Choose two additional points $X$ and $Y$ such that $X$ lies on the continuation of the segment $BZ$,  $Y$ lies on the continuation of the segment $AZ$, and $ \angle XAB=\angle YBA = \pi/2$. Consider a family of triangles  $\Omega_s=A BP_s$, $0 \le s\le 1$, such that 
$P_0 =X$, $P_1=Y$, and as $s$ changes, the point $P_s$ continuously moves along the union of segments $[X,Z]\cup[Z,Y]$. In particular, $X=P_0$, $Y=P_1$, and $Z=P_{s_0}$ for some $0<s_0<1$. We impose the Dirichlet condition on $[A,P_s]$ and the Neumann condition on $[B,P_s]$. The Steklov condition on $(A,B)$ does not change.

\begin{figure}[htb]
\begin{center}
\includegraphics{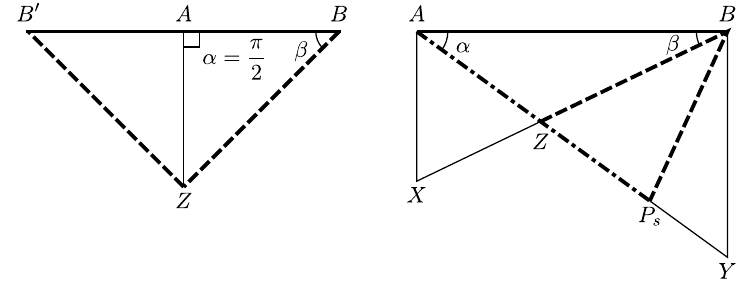}
\end{center}
\caption{Construction of auxiliary triangles in the case $\alpha=\pi/2$ (left) and  $\alpha<\pi/2$ (right)\label{fig:triangles}}
\end{figure}

For each of the triangles $\Omega_s$ one can construct a system of quasimodes similarly to the sloshing or Steklov-Dirichlet case; the only difference is that in this case we use the Dirichlet  Peters solutions in one angle and the Neumann Peters solutions in the other. 

It remains to show that this system is complete. In order to do that we use the same approach as in Lemma 
\ref{lem:monotonicityconsequences}. A simple modification of the argument for the domain monotonicity of the sloshing and Steklov-Dirichlet eigenvalues (see \cite{BKPS}) yields that 
for all $k\ge 1$, $\lambda_k(\Omega_s)\le \lambda_k(\Omega_{s'})$  if $s\le s'$. Indeed, if $P_s \in [X,Z]$ we can continue any trial function on a smaller triangle by zero across  $AP_s$ to get a trial function on a larger triangle, which means that the eigenvalues increase with $s$. Similarly, if $P_s\in [Z,Y]$, we can always restrict a trial function on a larger triangle to a smaller triangle, and this procedure decreases the Dirichlet energy without changing the denominator in the Rayleigh quotient. Hence, the eigenvalues increase with $s$ in this case as well. 

Since we  have already proved \eqref{mixedasymp} for right triangles, in the notation of  Lemma \ref{lem:monotonicityconsequences} it means that $k(j,0)=k(j,1)=j$ for all sufficiently large $j$. Following the same logic as in this  lemma we get that $k(j,s)=s$ for any $s$, in particular, for $s=s_0$. This completes the proof of the proposition. \qed


\section{Domains with curvilinear boundary}\label{sec:curvilinear}
\subsection{Domains which are straight near the surface} 
Suppose $\Omega$ is a sloshing domain whose sides are straight in a neighbourhood of $A$ and $B$. Let $T$ be the sloshing triangle with the same angles $\alpha$ and $\beta$ as $\Omega$; from the results in the previous Sections we understand sloshing eigenvalues and eigenfunctions on $T$ very well. As before, denote the sloshing eigenfunctions on $T$ by $\{u_j(z)\}$, 
with eigenvalues $\lambda_j$. Further, denote the original \emph{uncorrected} quasimodes on $T$, as in \eqref{eq:uncorrectedquasimodes}, by $\{v'_{\sigma_j}(z)\}$, with quasi-eigenvalues $\sigma_j$. Our strategy will be to use the quasimodes $v'_{\sigma_j}$ on $T$, cut off appropriately and modified slightly, as quasimodes for sloshing on $\Omega$.

Indeed, let $\chi(z)$ be a cut-off function on $\Omega$ with the following properties, which are easy to arrange:
\begin{itemize}
\item $\chi(z)$ is equal to $1$ in a neighbourhood of $S$;
\item $\chi(z)$ is supported on $\Omega\cap T$ and in particular is zero wherever the boundaries of $\Omega$ are not straight;
\item At every point $z\in\partial\Omega$, $\nabla\chi(z)$ is orthogonal to the normal $\boldsymbol{n}(z)$ to $\partial\Omega$.
\end{itemize}
Now define a set of functions on $\Omega$, by
\begin{equation}\label{eq:uncorrectedquasimodes2}
\tilde v'_{\sigma_j}(z):=\chi(z)v'_{\sigma_j}(z).
\end{equation}

The functions $\tilde v'_{\sigma_j}(z)$ do not quite satisfy Neumann conditions on the boundary, for the same reason the $v_{\sigma_j}(z)$ do not; we will have to correct them. Note, however, that since $\nabla\chi(z)\cdot n(z)=0$,
\[
\frac{\partial}{\partial n}\tilde v'_{\sigma_j}(z)=\chi(z)\frac{\partial}{\partial n}\tilde v'_{\sigma_j}(z),\qquad z\in\W,
\]
and hence, analogously to \eqref{eq:gradcorrect},
\[
\|\frac{\partial}{\partial n}\tilde v'_{\sigma_j}\|_{C^0(W)}\leq C\sigma_j^{-\mu}\quad\text{ on }\W.
\]
We correct these in precisely the same way as before, by adding a function $\tilde\eta_{\sigma_j}$ which is harmonic and which has Neumann data bounded everywhere in $C^0$ norm by $C\sigma_j^{-\mu}$. As before, since the Neumann-to-Dirichlet map is bounded on $L^2_*(S)$, we have
\[
\|\frac{\partial}{\partial y}\tilde\eta_{\sigma_j}\|_{L^2(S)}+\|\tilde\eta_{\sigma_j}\|_{L^2(S)}\leq C\sigma_j^{-\mu}.
\]
Then we define corrected quasimodes
\[
\tilde v_{\sigma_j}(z):=\tilde v'_{\sigma_j}(z)-\tilde \eta_{\sigma_j}(z),
\]
which now satisfy Neumann boundary conditions.

We may also compute their Laplacian. Since $v'_{\sigma_j}(z)$ are harmonic functions and $\tilde\eta_{\sigma_j}(z)$ are also harmonic, we have
\begin{equation}\label{eq:productruleformula}
\Delta\tilde v_{\sigma_j}(z)=\Delta\tilde v'_{\sigma_j}(z)=(\Delta\chi(z))v'_{\sigma_j}(z)+2\nabla\chi(z)\cdot\nabla v'_{\sigma_j}(z).
\end{equation}
Therefore $\Delta\tilde v_{\sigma_j}(z)$ is supported only where $\nabla\chi(z)$ is nonzero. For later use, we need to estimate the $L^2$ norm of $\Delta\tilde v_{\sigma_j}(z)$. It is immediate from the formula \eqref{eq:productruleformula} that for some universal constant $C$ depending on the cutoff function,
\[
\|\Delta\tilde v_{\sigma_j}(z)\|_{L^2(\Omega)}\leq C\|v'_{\sigma_j}(z)\|_{H^1(T\,\cap\,\mbox{supp}(\nabla\chi))}.
\]
The discussion before \eqref{eq:gradcorrect} shows that the $C^1$ norm of $v'_{\sigma_j}(z)$ is bounded by $C\sigma_j^{-\mu}$, and therefore obviously
\begin{equation}\label{eq:estimatesonetatilde}
\|v'_{\sigma_j}(z)\|_{H^1(T\,\cap\,\mbox{supp}(\nabla\chi))}\leq C\sigma_j^{-\mu}.
\end{equation}
Thus
\begin{equation}\label{eq:estimate77}
\|\Delta\tilde v_{\sigma_j}(z)\|_{L^2(\Omega)}\leq C\sigma_j^{-\mu}.
\end{equation}

We would like to use $\{\tilde v_{\sigma_j}|_S\}$ as our quasimodes. However, since $\tilde v_{\sigma_j}$ are not harmonic functions in the interior, we do not have $\mathcal D(\tilde v_{\sigma_j}|_S)=(\frac{\partial}{\partial_y}\tilde v_{\sigma_j})|_S$. We need to correct $\tilde v$ to be harmonic. To do this, we must solve the mixed boundary value problem
\begin{equation}\label{eq:adjunctproblem}
\begin{dcases} \Delta\tilde\varphi_{\sigma_j}=-\Delta\tilde v_{\sigma_j}&\quad\text{ on }\Omega;\\ 
\frac{\partial\tilde\varphi_{\sigma_j}}{\partial n}=0&\quad\text{ on }\W;\\ 
\tilde\varphi_{\sigma_j}=0&\quad\text{ on }S,
\end{dcases}
\end{equation}
as then we can use 
\[
\bar v_{\sigma_j}:=\tilde v_{\sigma_j}+\tilde\varphi_{\sigma_j}
\] 
as our quasimodes. The key estimate we want is given by
\begin{lemma}\label{lem:estimate78} If $\tilde\varphi_{\sigma_j}$ solves \eqref{eq:adjunctproblem}, then for some constant $C$ independent of $\sigma_j$,
\[
\left\|\frac{\partial}{\partial n}\tilde\varphi_{\sigma_j}\right\|_{L^2(S)}+\|\tilde\varphi_{\sigma_j}\|_{L^2(S)}\leq C\sigma^{-\mu}.
\]
\end{lemma}
\begin{proof}
The theory of mixed boundary value problems on domains with corners is required here; originally due to Kondratiev \cite{Kondr}, it has been developed nicely by Grisvard \cite{grisvard} in the setting of exact polygons. We will therefore transfer our problem to an exact polygon, namely $T$, via a conformal map. Let $\Phi:\Omega\to T$ be the conformal map which preserves the vertices and takes an arbitrary point $Z'$ on $\W$ to the vertex $Z$ of the triangle $\triangle ABZ$. Conformal maps preserve Dirichlet and Neumann boundary conditions, so the problem \eqref{eq:adjunctproblem} becomes, with $\varphi=\Phi_*(\tilde\varphi)$,
\begin{equation}\label{eq:conformaladjunctproblem}
\begin{dcases} \Delta\varphi_{\sigma_j}=-|\Phi'(z)|^2\Delta\tilde v_{\sigma_j}(z)&\quad\text{ on }T ;\\ 
\frac{\partial\varphi_{\sigma_j}}{\partial n}=0&\quad\text{ on }\W_T;\\
\varphi_{\sigma_j}=0&\quad\text{ on }S_T.
\end{dcases}
\end{equation}
By  \cite[Lemma 4.4.3.1]{grisvard}, there is indeed a unique solution $\varphi_{\sigma_j}(z)\in H^1(T)$ to \eqref{eq:conformaladjunctproblem}. Moreover, in this setting, the solution is actually in $H^2(T)$. Indeed, since all of the angles of $T$ are less than $\pi/2$ and in particular we do not have any corners of angle exactly $\pi/2$,   \cite[Theorem 4.4.3.13]{grisvard} applies. It tells us that $\varphi_{\sigma_j}(z)$ equals an element of $H^2(T)$ plus a linear combination of explicit separated-variables solutions at the corners, which Grisvard denotes $S_{j,m}$. However, consulting the definition of $S_{j,m}$ in  of \cite[equation (4.4.3.7)]{grisvard}, it is immediate that   each $S_{j,m}$ is itself in $H^2(T)$ (see also \cite[Lemma 4.4.3.5]{grisvard}) and therefore $\varphi_{\sigma_j}(z)\in H^2(T)$.

Now that \eqref{eq:conformaladjunctproblem} has a solution $\varphi_{\sigma_j}(z)\in H^2(T)$, we may apply the \emph{a priori} estimate  \cite[(4.1.2)]{grisvard}. It tells us that
\[
\|\varphi_{\sigma_j}\|_{H^2(T)}\leq C\left(\|\Delta\varphi_{\sigma_j}\|_{L^2(T)}+\|\varphi_{\sigma_j}\|_{L^2(T)}\right)\leq C\|\Delta\varphi_{\sigma_j}\|_{L^2(T)},
\]
where the second inequality follows from the fact that the first eigenvalue of  of the Laplacian on $T$ with Dirichlet conditions on $S_T$ and Neumann conditions on $\W_T$ is positive.
But since $\varphi_{\sigma_j}(z)$ solves \eqref{eq:conformaladjunctproblem},
\[\|\varphi_{\sigma_j}\|_{H^2(T)}\leq C\|\,|\Phi'(z)|\Delta\tilde v_{\sigma_j}(z)\|_{L^2(T)}.\]
Now observe that $|\Phi'(z)|$ is certainly smooth and bounded on the support of $\Delta\tilde v_{\sigma_j}(z)$, since that support is away from the vertices. Therefore, by \eqref{eq:estimate77},
\begin{equation}\label{eq:onemorebound}
\|\varphi_{\sigma_j}(z)\|_{H^2(T)}\leq C\sigma^{-\mu}.
\end{equation}

This can now be used to complete the proof. By \eqref{eq:onemorebound}, using the definition of Sobolev norms and the trace restriction Theorem,
\[
\left\|\frac{\partial}{\partial y}\varphi_{\sigma_j}\right\|_{L^2(S)}\leq 
\left\|\frac{\partial}{\partial y}\varphi_{\sigma_j}\right\|_{H^{1/2}(S)}
\leq C\left\|\frac{\partial}{\partial y}\varphi_{\sigma_j}\right\|_{H^1(T)}
\leq C\|\varphi_{\sigma_j}\|_{H^2(T)}\leq C\sigma^{-\mu}.
\]

All that remains to prove Lemma \ref{lem:estimate78} is to undo the conformal map. By \cite[Lemma 4.3]{PP2014}  (technically the local version, but the proof is local), since $\Phi$ preserves the angles at $A$ and $B$, $\Phi$ is $C^{1,\alpha}$ for any $\alpha$, and hence $C^{\infty}$, in a neighbourhood of $S$, and $\Phi'$ is nonzero. Note that $\Phi$ may not be $C^{\infty}$ near $Z'$ in particular, but we do not care about that region. Since $\Phi$ is smooth in a neighbourhood of $S$ and $\Phi'$ is nonzero, both the measure on $S$ and the magnitude of the normal derivative $\frac{\partial}{\partial y}$ change only up to a constant bounded above and below, and hence, as desired,
\[
\left\|\frac{\partial}{\partial n}\tilde\varphi_{\sigma_j}\right\|_{L^2(S)}\leq C\sigma_j^{-\mu}.
\]
An identical argument without taking $\frac{\partial}{\partial y}$ shows that
\[\|\tilde\varphi_{\sigma_j}\|_{L^2(S)}\leq C\sigma_j^{-\mu};\]
in fact, the bound is actually on the $H^{3/2}(S)$ norm. Hence the sum is bounded by $C\sigma_j^{-\mu}$ as well, which proves Lemma \ref{lem:estimate78}.
 \end{proof}

Now we claim our putative quasimodes $\bar v_{\sigma_j}$ on $\Omega$ satisfy a quasimode estimate:
\begin{lemma} There is a constant $C$ such that
\begin{equation}\label{eq:quasiestimate79}
\left\|\mathcal{D}_{\Omega}\left(\left.\bar v_{\sigma_j}\right|_S\right)-\sigma_j \bar v_{\sigma_j}\right\|_{L^2(S)}\leq C\sigma_j^{1-\mu}.
\end{equation}
\end{lemma}
Indeed, $\bar v_{\sigma_j}$  are harmonic and satisfy Neumann conditions on $\W$, so we know that 
\[
\mathcal{D}_{\Omega}\left(\left.\bar v_{\sigma_j}\right|_S\right)=\left.\frac{\partial}{\partial y}\bar v_{\sigma_j}\right|_S.
\]
As a result,
\[
\mathcal{D}_{\Omega}\left(\left.\bar v_{\sigma_j}\right|_S\right)-\sigma_j \bar v_{\sigma_j}
=\left.\frac{\partial}{\partial y}\tilde v_{\sigma_j}\right|_S-\sigma_j\tilde v_{\sigma_j}|_S+\left.\frac{\partial}{\partial y}\tilde\varphi_{\sigma_j}\right|_S-\sigma_j\tilde\varphi_{\sigma_j}|_S,
\]
which by definition of $\tilde v_{\sigma_j}$ is
\[
\left.\frac{\partial}{\partial y}\tilde v'_{\sigma_j}\right|_S
-\sigma_j\tilde v'_{\sigma_j}|_S
-\left.\frac{\partial}{\partial y}\tilde\eta_{\sigma_j}\right|_S
+\sigma_j\tilde\eta_{\sigma_j}|_S
+\left.\frac{\partial}{\partial y}\tilde\varphi_{\sigma_j}\right|_S-\sigma_j\tilde\varphi_{\sigma_j}|_S.
\]
The last two terms have $L^2$ norms which are bounded by $C\sigma^{1-\mu}$, by Lemma \ref{lem:estimate78}. The third and fourth terms satisfy the same bound as a consequence of \eqref{eq:estimatesonetatilde}. For the first two terms, recall that $\tilde v'_{\sigma_j}=\chi v'_{\sigma_j}$ and $\chi=1$ near $S$, so we need to estimate the $L^2$ norm of
\[
\left.\frac{\partial}{\partial y}v'_{\sigma_j}\right|_S-\sigma_jv'_{\sigma_j}|_S.
\]
However, $v'_{\sigma_j}$ are the explicit uncorrected quasimodes from \eqref{eq:uncorrectedquasimodes}. By the proof of Lemma \eqref{lem:keyquasi}, these two terms both have $L^2$ norms bounded by $C\sigma_j^{1-\mu}$, completing the proof of \eqref{eq:quasiestimate79}.

Now Theorem \ref{thm:linalg} applies immediately to tell us that near each sufficiently large quasi-frequency $\sigma_j$, there exists at least one sloshing eigenvalue 
$\lambda_k(\Omega)$. The problem, again, is to prove completeness.

\subsection{Completeness} We use the same strategy that we did for our original quasimode problem. Since the strategy is the same, we omit some details.

First we prove completeness when $\alpha=\beta=\pi/2q$. To do this, we instead use Hanson-Lewy quasimodes. All the analysis we have done in the previous subsection goes through with $C\er^{-\sigma/C}$ replacing $C\sigma^{-\mu}$ everywhere, and we get \eqref{eq:quasiestimate79} with a right-hand side of $C\er^{-j/C}$. So Theorem \ref{thm:linalg} applies. It shows in particular that there are (linear combinations of) sloshing eigenfunctions on $\Omega$, which we call $w_{j,\Omega}$, with eigenvalues exponentially close to $\sigma_j$ and with the property that
\[
\|w_{j,\Omega}-\bar v_{\sigma_j}\|_{L^2}\leq C\er^{-j/C}\mbox{ as }j\to\infty.
\]
However $\|\bar v_{\sigma_j}-\tilde v_{\sigma_j}\|_{L^2(S)}$, $\|\tilde v_{\sigma_j}-v'_{\sigma_j}\|_{L^2(S)}$, $\|v'_{\sigma_j}-v_{\sigma_j}\|_{L^2(S)}$, and $\|v_{\sigma_j}-u_j\|_{L^2(S)}$ are all bounded by $C\er^{-j/C}$ themselves in this case, so we have in all that
\[
\|w_{j,\Omega}-u_{j}\|_{L^2}\leq C\er^{-j/C}\mbox{ as }j\to\infty.
\]

We know that $\{u_j\}$ are complete. By the same argument as in Lemma \ref{lem:barykrein}, the sequence $\{w_{j,\Omega}\}$ must at some point contain only pure eigenfunctions, and for suitably large $N$ the sequences $\{w_{j,\Omega}\}_{j\geq N}$ and $\{u_{j}\}_{j\geq N}$ are missing the same number of eigenfunctions, which shows that
\[\lambda_{j,\Omega}=\lambda_{j,T}+O(\er^{-j/C}).\]
This proves completeness, along with full exponentially accurate eigenvalue asymptotics, in the $\alpha=\beta=\pi/2q$ case.

Finally, we use domain monotonicity to deal with general $\alpha$ and $\beta$; the argument runs exactly as in subsection \ref{subsection:monotonicity}. The result of Davis \cite{dav65} shows completeness in the $\pi/2$ case, and we can sandwich our domain $\Omega$ as $\Omega_{1/2}$ in a continuous, nested family of sloshing domains $\{\Omega_s\}$ that satisfy the hypotheses of Lemma \ref{lem:monotonicityconsequences}. Note that Lemma \ref{lem:monotonicityconsequences} is stated in such a way as to apply here as well.

\subsection{General curvilinear domains}
\label{subs:curvilin}
We now generalise further, to domains with curvilinear boundary. The authors are grateful to Lev Buhovski for suggesting the approach in this subsection \cite{Buh}.

Suppose that $\Omega$ is a sloshing domain which satisfies the following conditions:
\begin{enumerate}
\item[(C1)] $\Omega$ is simply connected;
\item[(C2)] The walls $\W$  are Lipschitz;
\item[(C3)] For any $\epsilon>0$, there exist sloshing domains $\Omega_{-}$ and $\Omega_+$ with piecewise smooth boundary, with $\Omega_{-}\subset\Omega\subset\Omega_+$, and with the additional property that $\Omega_{-}$ and $\Omega_+$ are straight lines in a $\delta-$neighbourhood of $A$ and $B$ with vertex angles in $[\alpha-\epsilon,\min\{\alpha+\epsilon,\pi/2\}]$.
\end{enumerate}
The point of this third condition is that we have already proved sloshing eigenvalue asymptotics for $\Omega_-$ and $\Omega_+$, which we will need in the subsequent domain monotonicity argument. This third condition is satisfied, for example, if the boundary $S$ is $C^1$ in any small neighbourhood of the vertices $A$ and $B$ and the angles are strictly less than $\pi/2$. It is also satisfied under the ``local John's condition" (see Propositions \ref{prop:pi2} and \ref{prop:pi22}) if one or both angles equals $\pi/2$, and makes clear why that condition is necessary: we have not proved sloshing asymptotics in the case where one or both angles are greater than $\pi/2$.

Under these conditions, we claim the asymptotics \eqref{eq:sloshingasymptotics}, \eqref{eq:sloshingasymptoticsdirichlet}, as well as \eqref{pi2} and \eqref{pi22} for domains satisfying the additional ``local John's condition', and thereby complete the proofs of Theorems \ref{thm:main} and \ref{thm:Dirichlet} as well as Propositions \ref{prop:pi2} and \ref{prop:pi22}.

\begin{proof} Assume $L=1$ for simplicity. Pick any $\gamma>0$. By continuity, there exists a sufficiently small $\epsilon$ such that for any $\epsilon'\in[-\epsilon,\epsilon]$,
\[
\left|\frac{1}{\alpha+\epsilon'}-\frac{1}{\alpha}\right|+\left|\frac{1}{\beta+\epsilon'}-\frac{1}{\beta}\right|<\frac{8\gamma}{2\pi^2}.
\]

By condition (C3) above, there exist $\Omega_-\subset\Omega$ and $\Omega_+\supset\Omega$ as described. By the definition of $\epsilon$, and the explicit form of $\sigma_j$,
\[
|\sigma_j(\Omega_+)-\sigma_j(\Omega)|<\gamma/2\mbox{ and }|\sigma_j(\Omega_-)-\sigma_j(\Omega)|<\gamma/2\text{ for all }j.
\]
But by our previous work, we have eigenvalue asymptotics for $\Omega_+$ and $\Omega_-$. In particular, there exists $N$ so that for all $j\geq N$, 
\[
|\lambda_j(\Omega_+)-\sigma_j(\Omega_+)|<\gamma/2\text{ and }|\lambda_j(\Omega_-)-\sigma_j(\Omega_-)|<\gamma/2.
\]
And by domain monotonicity (for Neumann --- inequalities reverse for the Dirichlet case),
\[
\lambda_j(\Omega_-)\leq\lambda_j(\Omega)\leq\lambda_j(\Omega_+)\text{ for all }j.
\]
Putting these ingredients together, for $j\geq N$,
\[
\lambda_j(\Omega)\leq\lambda_j(\Omega_+)<\sigma_j(\Omega_+)+\frac{\gamma}{2}\leq\sigma_j(\Omega)+\gamma.
\]
Similarly,
\[
\lambda_j(\Omega)\geq\lambda_j(\Omega_-)>\sigma_j(\Omega_-)-\frac{\gamma}{2}\geq\sigma_j(\Omega)-\gamma.
\]
We conclude that for $j\geq N$, $|\lambda_j(\Omega)-\sigma_j(\Omega)|<\gamma$. Since $\gamma>0$ was arbitrary, we have the asymptotics with $o(1)$ error that we want.
\end{proof}



\appendix

\section{Proof of Theorem \ref{wedgemodels}}\label{proofpeters}
In this section, we prove Theorem \ref{wedgemodels}. The proof follows \cite{peters} for the most part, doing the extra work needed to prove the careful remainder estimates. There is a key difference of notation: throughout, we use $\mu_\alpha=\mu=\pi/(2\alpha)$, where Peters uses $\mu=\pi/\alpha$.
\subsection{Robin-Neumann problem}
Following Peters, we complexify our problem by setting $z=\rho \er^{\ir \theta}$ and look for an analytic function $f(z)$ in $\Sct{\alpha}$ with $\phi=\Re (f)$. The boundary conditions must be rewritten in terms of $f$. Using the Cauchy-Riemann equations we have $\phi_y=\Re (\ir f'(z))$ and $\phi_x=\Re (f'(z))$. After some algebraic transformations, which we skip for brevity, the problem \eqref{modelproblem:neumann} becomes
\begin{equation}\label{analyticmodel:neumann}
\begin{dcases}f(z)&\quad\text{ is analytic in }\Sct{\alpha},\\ 
\Re (\ir f'(z))=\Re(f(z))&\quad\text{ for }z\in\Sct{\alpha} \cap\{\theta=0\},\\
\Re (\ir\er^{-\ir \alpha}f'(z))=0&\quad\text{ for }z\in \Sct{\alpha}\cap\{\theta=-\alpha\}.\\
\end{dcases}
\end{equation}

\subsection{Peters solution}

We now write down Peters solution, in a form due to Alker \cite{Alker}. Note that Alker's $y$-axis points in the opposite direction from Peters and so we have modified the expression in \cite{Alker} accordingly. First define the auxiliary function
\begin{equation}\label{eq:Iintegral}
I_{\alpha}(\zeta):=\frac{1}{\pi}\int_0^{\zeta\infty}\log\big(1+v^{-\pi/\alpha}\big)\frac{\zeta}{v^2+\zeta^2}\,\dr v.
\end{equation}
Note that $I_{\alpha}(\zeta)$ is defined for $\arg(\zeta)\in(-\alpha,\alpha)$, which includes the real axis. As in the appendix of \cite{peters}, $I_{\alpha}(\zeta)$ has a meromorphic continuation, with finitely many branch points, to the entire complex plane, which we also call $I_{\alpha}(\zeta)$. Each branch point is logarithmic; there is one at the origin, and others on the unit circle in the negative real half-plane. We then let
\begin{equation}\label{eq:gorig}
g_{\alpha}(\zeta)=\exp\left(-I_{\alpha}(\zeta \er^{-\ir \alpha})+\log\left(\frac{\zeta+\ir}{\zeta}\right)\right).
\end{equation}
The function $g_{\alpha}(\zeta)$ is originally defined for $\arg(\zeta)\in(0,2\alpha)$, but has a meromorphic continuation to the entire complex plane minus a single branch cut from the origin, with singularities along the portion of the unit circle outside the sector $-\pi/2-\alpha\leq\arg(\zeta)\leq\pi/2+\alpha$. For $\Re(\zeta)>0$ we have the representation formula \cite[(4.8)]{peters}
\begin{equation}\label{eq:g1}
g_{\alpha}(\zeta)=\exp\left(-\frac{1}{\pi}\int_0^{\infty}\log\left(\frac{1-t^{-2\mu}}{1-t^{-2}}\right)\frac{\zeta}{t^2+\zeta^2}\ \dr t\right).
\end{equation}
Finally, let $P$ be a keyhole path, consisting of the union of a nearly full circle of radius 2, traversed counterclockwise, and two linear paths, one on each side of the angle $\theta=\pi+\alpha/2$. 
Then the Peters solution is given by
\begin{equation}\label{eq:solution}
f(z)=\frac{\mu^{1/2}}{\ir\pi}\int_{P}\frac{g_\alpha(\zeta)}{\zeta+\ir}\er^{z\zeta}\,\dr\zeta.
\end{equation}

\subsection{Verification of Peters solution}

Let us verify that Peters solution, which is obviously analytic, is actually a solution by checking the boundary conditions in \eqref{analyticmodel:neumann}, beginning with the one along the real axis. This follows \cite[section 6]{peters}.

By \cite[p. 335]{peters}, $g_{\alpha}(\zeta)\to 1$ as $\zeta\to\infty$ (see Proposition \ref{prop:asymptoticsofg} for a rigorous proof), and since $\er^{z\zeta}$ is very small as $\zeta\to\infty$ along $P$, differentiation under the integral sign in \eqref{eq:solution} is justified. We get
\begin{equation}\label{eq:diff1}
\ir f'(z)-f(z)=\frac{\mu^{1/2}}{\ir\pi}\int_{P}\frac{g_{\alpha}(\zeta)(\ir\zeta-1)}{\zeta+\ir}\er^{z\zeta}\,\dr\zeta=\frac{\mu^{1/2}}{\pi}\int_{P}g_{\alpha}(\zeta)\er^{z\zeta}\,\dr\zeta.
\end{equation}
We claim this integral is pure imaginary whenever $z$ is on the positive real axis. Indeed, we may shift the branch cut of $g_{\alpha}(\zeta)$ so that it lies along the negative real axis, and shift $P$ so that it is symmetric with respect to the real axis. Observe that $g_{\alpha}(\zeta)$ is real whenever $\zeta$ is real and positive. By reflection, this implies that $g_{\alpha}(\overline\zeta)=\overline{g_{\alpha}(\zeta)}$ for all $\zeta$. The analogous statement is true for $\er^{z\zeta}$ and hence for $g_{\alpha}(\zeta)\er^{z\zeta}$, and it follows immediately from symmetry of $P$ that \eqref{eq:diff1} is purely imaginary. Thus $\Re(\ir f'(z)-f(z))=0$, as desired.

For the other boundary condition, we compute
\begin{equation}\label{eq:diff2}
\ir\er^{-\ir \alpha}f'(z)=\frac{\mu^{1/2}}{\pi}\int_{P}\frac{g_{\alpha}(\zeta)\er^{-\ir \alpha}\zeta}{\zeta+\ir}\er^{z\zeta}\,\dr \zeta,
\end{equation}
and claim that this is pure imaginary when $\arg(z)=-\alpha$ - i.e. when $z=\rho \er^{-\ir \alpha}$ for some $\rho>0$. Shifting the branch cut and contour of integration, and letting $w=\zeta \er^{-\ir \alpha}$, the integral \eqref{eq:diff2} becomes
\begin{equation}\label{eq:diff3}
\frac{\mu^{1/2}}{\pi}\int_{P_1}\frac{g_{\alpha}(w \er^{\ir \alpha})w \er^{\ir \alpha}}{w \er^{\ir \alpha}+\ir}\er^{\rho w}\,\dr w=
\frac{\mu^{1/2}}{\pi}\int_{P_1}\exp(-I_{\alpha}(w))\er^{\rho w}\,\dr w,
\end{equation}
where we have chosen $P_1$ so that it is symmetric with respect to the real axis and the branch cut is along the negative real axis. The boundary condition follows immediately as above, since from \eqref{eq:Iintegral} $I_{\alpha}(w)$ is positive for $w$ real and positive.

\subsection{Asymptotics of Peters solution as $z\to\infty$}

Now we study the asymptotics as $z\to\infty$ of \eqref{eq:solution}, following \cite[section 7]{peters} but with more rigour. Let, as before,
\[
\chi=\chi_{\alpha,N}=\frac{\pi}{4}(1-\mu)=\frac{\pi}{4}\left(1-\frac{\pi}{2\alpha}\right).
\]
\begin{theorem}\label{Peters} There exist a function $E:(0,\pi/2)\to\mathbb R$ and a complex-valued function $R_{\alpha}(z)$ depending on $\alpha$ such that
\begin{equation}
f(z)=\er^{E(\alpha)}\er^{-\ir (z-\chi)}+R_{\alpha}(z),
\end{equation}
where for any fixed $\alpha\in (0,\pi/2)$ there exists a constant $C$ such that for all $z\in \Sct{\alpha}$,
\begin{equation}
|R_{\alpha}(z)|\leq Cz^{-\mu},\ |\nabla_z R_{\alpha}(z)|\leq Cz^{-\mu-1}.
\end{equation}
\end{theorem}
\begin{remark} A slightly stronger version of this theorem (identifying $E(\alpha)$) is claimed in \cite{Alker}. However, the proof is not given and the extraction does not seem obvious, although it is numerically clear. So we prove this version, which is all we need since we may scale by an overall constant anyway. \end{remark}
\begin{remark} For the solution $\phi=\Re(f(z))$ to our original problem, we see that
\begin{equation}
\phi(x,0)=\er^{E(\alpha)}\cos(z-\chi)+R_{\alpha}(x+0\ir),
\end{equation}
which is the radiation condition we wanted in the first place.
\end{remark}

\begin{proof} 
 The portion of the integral \eqref{eq:solution} along the infinite line segments decays exponentially in $|z|$ (note that $\arg(z\zeta)$ is between $\pi-\alpha/2$ and $\pi+\alpha/2$). So we may deform the contour $P$ to a contour $P'$ consisting of the union of a circle of radius $1/2$ and line segments along each side of the branch cut. This deformation passes through singularities of $g_{\alpha}(\zeta)/(\zeta+\ir)$, one at $\zeta=-\ir$ and others at $\zeta=\lambda_1,\ldots,\lambda_m$ along the unit circle in the negative real half-plane, with $\arg(\lambda_j)\notin(-\pi/2-\alpha,\pi/2+\alpha)$. We see that
\begin{equation}\label{eq:solutionmark1}
\begin{split}
f(z)&=2\mu^{1/2}g_{\alpha}(-\ir)\er^{-\ir z}+\frac{\mu^{1/2}}{\ir\pi}\int_{P'}\frac{g_{\alpha}(\zeta)}{\zeta+\ir}\er^{z\zeta}\,\dr \zeta\\
&\quad+\sum_{j=1}^m\er^{z\lambda_m}2\mu^{1/2}\mathop{\operatorname{Res}}\limits_{\zeta=\lambda_m}(g_{\alpha}(\zeta)/(\zeta+\ir)).
\end{split}
\end{equation}
Let
\begin{equation}\label{eq:solutionmark1a}
\begin{split}
\widetilde R_{\alpha}(z)&:=\frac{\mu^{1/2}}{\ir\pi}\int_{P'}\frac{g_{\alpha}(\zeta)}{\zeta+\ir}\er^{z\zeta}\,\dr \zeta;\\ 
R_{\alpha}(z)&:=\widetilde R_{\alpha}(z)+\sum_{j=1}^m\er^{z\lambda_m}2\mu^{1/2}\mathop{\operatorname{Res}}\limits_{\zeta=\lambda_m}(g_{\alpha}(\zeta)/(\zeta+\ir)).
\end{split}
\end{equation}
To prove Theorem \ref{Peters} we must first prove the error estimates for $R_{\alpha}(z)$ and then compute $g_{\alpha}(-\ir)$.

\subsubsection{Error estimates for $R_{\alpha}(z)$} Consider $R_{\alpha}(z)$. Since $\arg(z)\in(-\alpha,0)$, the finite sum of exponentials and its gradient decay exponentially as $|z|\to\infty$. The decay is actually uniform in $\alpha$ for $\alpha$ in any interval $(\epsilon,\pi-\epsilon)$ with $\epsilon>0$ fixed, as the number of residues is bounded in such an interval as well. Only $\widetilde R_{\alpha}(z)$ remains. Using \eqref{eq:Iintegral}, we may write
\begin{equation}\label{firstterm}
\widetilde R_{\alpha}(z)=\frac{\mu^{1/2}}{\ir\pi}\int_{P'}\er^{z\zeta}\exp(-I_{\alpha}(\zeta \er^{-\ir \alpha}))\frac{\dr \zeta}{\zeta}.
\end{equation}
The only singularity inside $P'$ is along the branch cut, including the branch point $0$. We see that we must understand the asymptotics of $I_{\alpha}(\zeta)$ as $\zeta\to 0$. Peters identifies the leading order term of $I_{\alpha}(\zeta)$ at $\zeta=0$ as $-\mu\log\zeta$, but we need to understand the remainder and prove bounds on it and its gradient. Throughout, we choose to suppress the $\alpha$ subscripts. The key is the following lemma:

\begin{lemma}\label{lem:A4} 
For $\zeta$ in a small neighbourhood of zero, with argument up to the branch cut on either side, we have
\[I(\zeta)=-\mu\log\zeta+p(\zeta)+R(\zeta),\]
where $p(\zeta)$ is a polynomial in $\zeta$ with $p(0)=0$, and where there is a constant $C$ such that
\[|R(\zeta)|\leq C|\zeta|^{2\mu}.\]
\end{lemma}
\begin{proof} The proof proceeds in two steps. First we prove this for all $\zeta$ with argument in a compact subset of $(-\pi/2,\pi/2)$. Then we use a functional relation satisfied by $I(\zeta)$ to extend to $\zeta$ with arguments up to and into the negative half-plane.

For the first step, we use the representation for $I(\zeta)$ on p. 353 of \cite{peters}. From the last paragraph of this page, we have
\[I(\zeta)=-\mu\log\zeta+\int_0^{1/2}\ln(1+v^{2\mu})\frac{\zeta}{v^2+\zeta^2}\,\dr v+\int_{1/2}^{\infty}\ln(1+v^{2\mu})\frac{\zeta}{v^2+\zeta^2}\,\dr v.\]
The last term is convergent, as $\ln(1+v^{2\mu})\sim 2\mu\ln v$ as $v\to\infty$, and differentiation in $\zeta$ under the integral sign is easily justified. So the last term is holomorphic, and by direct substitution it is zero at $\zeta=0$. For the second term, we use a Taylor expansion for $\ln(1+v^{2\mu})$, which is convergent when $v<1$, and obtain
\[\int_0^{1/2}\sum_{n=1}^{\infty}(-1)^{n+1}\frac{v^{2\mu n}}{n}\frac{\zeta}{v^2+\zeta^2}\,\dr v.\]
Change variables in the integral to $w=\zeta v$; we end up with
\begin{equation}\label{eq:changedvarint}\int_0^{1/2\zeta}\sum_{n=1}^{\infty}(-1)^{n+1}\frac{1}{n}\zeta^{2\mu n}w^{2\mu n}\frac{1}{w^2+1}\,\dr w.\end{equation}

Now further break up the integral \eqref{eq:changedvarint}, at $w=2$. The integral from $0$ to $2$ is bounded in absolute value by
\[
\int_0^2\sum_{n=1}^{\infty}|\zeta|^{2\mu n}\frac{2^{2\mu n}}{n}\,\dr w=2\sum_{n=1}^{\infty}|\zeta|^{2\mu n}\frac{2^{2\mu n}}{n}\leq 2\sum_{n=1}^{\infty}|2\zeta|^{2\mu n}=\frac{|4\zeta|^{2\mu}}{1-|2\zeta|^{2\mu}},
\]
which is bounded by $C|\zeta|^{2\mu}$ for some universal constant $C$ and sufficiently small $|\zeta|$. For the remainder of the integral, we use the Taylor expansion of $(w^2+1)^{-1}$ about infinity, which is valid for $w>1$: it is $w^{-2}-w^{-4}+w^{-6}+\dots$, and we end up with
\[
\int_2^{1/2\zeta}\sum_{n=1}^{\infty}(-1)^{n+1}\frac{1}{n}\zeta^{2\mu n}w^{2\mu n}\left(\sum_{m=1}^{\infty}(-1)^{m+1}w^{-2m}\right)\,\dr w.
\]
Now observe that $|\zeta w|\leq\frac 12$ on this interval and $w>2$, so both the sums in $n$ and $m$ are absolutely convergent and in fact the double sum is absolutely convergent. This justifies all rearrangements as well as term-by-term integration. Our last remaining piece thus becomes
\[\sum_{n=1}^{\infty}\sum_{m=1}^{\infty}\frac{(-1)^{n+m}}{n}\zeta^{2\mu n}\int_2^{1/2\zeta}w^{2\mu n-2m}\,\dr w,\]
which equals
\[
\begin{split}
&\qquad\sum_{n=1}^{\infty}\sum_{m=1}^{\infty}\frac{(-1)^{n+1+m}}{n}\zeta^{2\mu n}\left.\left(\frac{w^{2\mu n-2m+1}}{2\mu n-2m+1}\right)\right|^{1/2\zeta}_2\\
&=\sum_{n=1}^{\infty}\sum_{m=1}^{\infty}\frac{(-1)^{n+1+m}}{n(2\mu n-2m+1)}\zeta^{2\mu n}\left(\zeta^{-2\mu n+2m-1}2^{-2\mu n+2m-1}-2^{2\mu n-2m+1}\right)\\
&=\sum_{n=1}^{\infty}\sum_{m=1}^{\infty}\frac{(-1)^{n+1+m}}{n(2\mu n-2m+1)}\left(\zeta^{2m-1}2^{-2\mu n+2m-1}-\zeta^{2\mu n}2^{2\mu n-2m+1}\right).
\end{split}
\]
The sum now has two terms. The second term, with a $\zeta^{2\mu n}$, can be summed first in $m$ (obviously) and then in $n$, and the whole thing is bounded by a multiple of the $n=1$ term, namely $C|\zeta|^{2\mu}$. The first term can be summed first in $n$, at which point it becomes an expression of the form $\sum_{m=1}^{\infty}a_m\zeta^{2m-1}$. Since the $a_m$ do not grow too fast as $m\to\infty$ --- in fact they grow as $2^{2m}$ --- this expression represents a holomorphic function in a disk of sufficiently small radius about the origin. Moreover this holomorphic function is zero at zero.

We have now shown that $I(\zeta)$ is the sum of $-\mu\log\zeta$, a holomorphic function zero at the origin, and a term bounded in absolute value by $C|\zeta|^{2\mu}$. Taking the Taylor series of the holomorphic function, separating out the finitely many terms which are not $O(|\zeta|^{2\mu})$, and calling them $p(\zeta)$ completes the proof.

It remains to extend the argument outside the positive real half-plane. Let $h(\zeta)=I(\zeta)+\mu\log\zeta$. Then we know that $h(\zeta)$ is holomorphic in a disk with the exception of a branch cut, and it continues across the branch cut because $I(\zeta)$ does. We also know, from the first line on p. 353 of \cite{peters}, that for all $\zeta$ away from the branch cut\label{page:foot}\footnote{The authors would like to thank Marcello Malagutti for pointing out that the exponential in the right-hand side of the formula was missing in the published version of the paper. This does not affect the statement of Lemma  \ref{lem:A4}, and its proof is adjusted accordingly.},
\[\exp(I(\zeta \er^{2\ir\alpha}))=\exp(I(\zeta))\frac{\zeta \er^{\ir \alpha}+\ir}{\zeta \er^{\ir \alpha}-\ir}.\]
Applying this gives an equivalent relation for $h(\zeta)$, noting that $\er^{2\ir\alpha\mu}=\er^{\ir \pi}=-1$:
\[\exp(h(\zeta \er^{2\ir\alpha}))=-\exp(h(\zeta))\frac{\zeta \er^{\ir \alpha}+\ir}{\zeta \er^{\ir \alpha}-\ir}.\]
The function $H(\zeta):=-(\zeta \er^{\ir \alpha}+\ir)/(\zeta \er^{\ir \alpha}-\ir)$ is holomorphic and nonzero in the disk, and we have for all nonzero $\zeta$:
\[\exp(h(\zeta \er^{2\ir\alpha}))=\exp(h(\zeta))H(\zeta).\]
Taking logarithms, we see that there is a branch $\log H(\zeta)$ for which
\[h(\zeta \er^{2\ir\alpha})=h(\zeta)+\log H(\zeta).\]

Recall that in the positive real half-plane (in a sector bounded away from the real axis),
\[h(\zeta)=p(\zeta)+R(\zeta).\]
Plugging in the relation for $h(\zeta)$, we see that for $-\pi/2+\epsilon<\arg\zeta<\pi/2-\epsilon$,
\[h(\zeta \er^{2\ir\alpha})=p(\zeta)+R(\zeta)+\log H(\zeta).\]
However, we are assuming that $\alpha<\pi/2$. Therefore, there is a range of $\zeta$, namely $-\pi/2+\epsilon<\arg\zeta<\pi/2-\epsilon-2\alpha$, where $\zeta \er^{2\ir\alpha}$ is in the positive half-plane and thus we have a second representation:
\[h(\zeta \er^{2\ir\alpha})=p(\zeta \er^{2\ir\alpha})+R(\zeta \er^{2\ir\alpha}).\]
Setting the previous two equations equal, we see that for this small range of $\zeta$:
\[p(\zeta \er^{2\ir\alpha})-p(\zeta)-\log H(\zeta)=R(\zeta)-R(\zeta \er^{2\ir\alpha}).\]
The left-hand side is holomorphic in a disk. The right-hand side is $O(|\zeta|^{2\mu})$. We conclude that the left-hand side must have a zero of order at least $2\mu$ at $\zeta=0$. Therefore, there is a $C$ such that for \emph{all} $\zeta$ in the disk,
\[p(\zeta \er^{2\ir\alpha})-p(\zeta)-\log H(\zeta)=O(|\zeta|^{2\mu}).\]
Rearranging, we observe that for $-\pi/2+\epsilon<\arg\zeta<\pi/2-\epsilon$,
\[h(\zeta \er^{2\ir\alpha}) = p(\zeta \er^{2\ir\alpha})+R(\zeta) - (p(\zeta \er^{2\ir\alpha})-p(\zeta)-\log H(\zeta)).\]
But we know both of the last two terms on the right-hand side are $O(\zeta^{2\mu})$ and thus
\[h(\zeta \er^{2\ir\alpha})=p(\zeta \er^{2\ir\alpha}) + O(\zeta^{2\mu})\]
as desired. This proves the lemma for $\zeta$ with arguments now up to $\pi/2-\epsilon+2\alpha$. Continuing this procedure, one step of size $2\alpha$ at a time, gives the result for all $\zeta$ with arguments in a neighbourhood of $[-\pi,\pi]$, completing the proof.
\end{proof}

Now we return to estimating the remainder. Using the Lemma and recalling $h(\zeta)=I(\zeta)+\mu\log\zeta$, we have
\begin{equation}\label{firstterm1a}
\widetilde R_{\alpha}(z)=-\frac{\mu^{1/2}}{\pi}\int_{P'}\er^{z\zeta}\zeta^{\mu}\er^{-h(\zeta)}\frac{d\zeta}{\zeta}.
\end{equation}
Change variables to let $w=z\zeta$; for $|z|>1$, the contour deforms smoothly back to $\arg(z)P'$ (which is $P'$ rotated by $\arg(z)$), as there are no singularities remaining inside $P'$. We get
\begin{equation}\label{firstterm2}
\widetilde R_{\alpha}(z)=-\frac{\mu^{1/2}z^{-\mu}}{\pi}\int_{\arg(z)P'}\er^ww^{\mu-1}\er^{-h(w/z)}\,\dr w.
\end{equation}
The contour can be deformed to two straight lines, one on either side of the branch cut. Since $h(\cdot)$ grows at most logarithmically along these lines (note that $\arg(w/z)=\pm\pi+\alpha/2$ on $\arg(z)P'$), the integral is bounded as $z\to\infty$. In fact, since $h(0)=0$, it converges to the corresponding integral with $\er^{-h(w/z)}$ replaced by $1$. Thus $|\widetilde R_{\alpha}(z)|\leq C z^{-\mu}$ as desired.

We must also estimate $|\nabla_z \widetilde R_{\alpha}(z)|$. However, $\widetilde R_{\alpha}(z)$ is holomorphic, so we just need to consider $\partial_z\widetilde R_{\alpha}(z)$. The differentiation brings down a factor of $\zeta$, which becomes $w/z$. The extra power of $w$ is absorbed, and the extra power of $z$ moves outside the integral, yielding decay of the form $z^{-\mu-1}$ instead of $z^{-\mu}$. The same analysis works for higher order derivatives, and this proves the remainder estimation part of Theorem \ref{Peters}. It remains only to evaluate $g_{\alpha}(-\ir)$.

\subsubsection{Evaluation of $g_{\alpha}(-\ir)$} We cannot simply use the formula \eqref{eq:g1} to evaluate $g_{\alpha}(-\ir)$, because that formula is only valid for $\zeta$ in the positive real half-plane. We must instead use the fact that
\begin{equation}\label{eq:backtog1}
\log g_{\alpha}(\zeta)=\log\left(\frac{\zeta+\ir}{\zeta}\right)-I(\zeta \er^{-\ir \alpha}),
\end{equation}
so we are interested in understanding the behaviour of $I(\zeta)$ near $\zeta=-\ir\er^{-\ir \alpha}$.

First we need a good representation for $I(\zeta)$ in this region. Considering the last equation on p. 350 of \cite{peters}, we see that for $\zeta$ with $-\pi/2-2\alpha<\arg(\zeta)<\pi/2+2\alpha$,
\begin{equation*}
\begin{split}
I(\zeta \er^{-\ir \alpha})&=\frac{1}{2\pi\ir}\left(\int_{I_1}\log(1+(-\ir u)^{-2\mu})\frac{\dr u}{u-\zeta \er^{-\ir \alpha}}+\int_{M_2}\log(1+(\ir u)^{-2\mu})\frac{\dr u}{u-\zeta \er^{-\ir \alpha}}\right)\\
&+\log\frac{\zeta+\ir}{\zeta}.
\end{split}
\end{equation*}
Here $I_1$ may be chosen to be the incoming path along $\arg(\zeta)=\pi/2$ and $M_2$ the outgoing path along $\arg(\zeta)=-\pi/2-2\alpha$ (each can be moved by up to $\alpha$ in either direction before running into a singularity of the integrand, but these are the most convenient choices). We would prefer that $I_1$ and $M_2$ be equal and opposite paths, but they are not. So we shift the path $I_1$ to the path $L_0$ which is incoming along the ray $\arg(\zeta)=\pi/2-2\alpha$. This is done for other reasons on p. 351 of \cite{peters}, and we pick up a contribution from the branch cut of the integrand. Overall we get
\begin{equation}\label{largerarearep}
\begin{split}
I(\zeta \er^{-\ir \alpha})&=\frac{1}{2\pi\ir}\left(\int_{L_0}\log(1+(-\ir u)^{-2\mu})\frac{\dr u}{u-\zeta \er^{-\ir \alpha}}+\int_{M_2}\log(1+(\ir u)^{-2\mu})\frac{\dr u}{u-\zeta \er^{-\ir \alpha}}\right)\\&+\log\frac{\zeta+\ir}{\zeta}-\log\frac{\zeta-\ir}{\zeta}.
\end{split}
\end{equation}
Using \eqref{eq:backtog1}, then plugging in $\zeta=-\ir$, shows that
\begin{equation}
\begin{split}
-\log(g_{\alpha}(-\ir))&=\frac{1}{2\pi\ir}\left(\int_{L_0}\log(1+(-\ir u)^{-2\mu})\frac{\dr u}{u+\ir \er^{-\ir \alpha}}+\int_{M_2}\log(1+(\ir u)^{-2\mu})\frac{\dr u}{u+\ir \er^{-\ir \alpha}}\right)\\
&-\log 2.
\end{split}
\end{equation}
Now we parametrise these integrals. For the first, let $u=t\ir\er^{-2i\alpha}$, and for the second, let $u=-t\ir\er^{-2i\alpha}$. With the appropriate signs, we find that
\begin{equation}\label{eq:evilintegral}
\begin{split}
-\log(g_{\alpha}(-\ir))&=-\log 2+\frac{1}{2\pi\ir}\int_0^{\infty}\log(1+t^{-2\mu})\left(\frac{1}{t-\er^{\ir \alpha}}-\frac{1}{t+\er^{\ir \alpha}}\right)\,\dr t
\\&=-\log 2+\frac{1}{\pi\ir}\int_0^{\infty}\log(1+t^{-2\mu})\frac{\er^{\ir \pi/2\mu}}{t^2-\er^{\ir \pi/\mu}}\,\dr t.
\end{split}
\end{equation}

We have now written $g_{\alpha}(-\ir)$ in terms of a convergent integral, that on the right of \eqref{eq:evilintegral}. The usual Taylor series expansions show that this integral is smooth in $\alpha$ for $\alpha\in(\epsilon,\pi-\epsilon)$. This integral seems nontrivial to evaluate and is not obviously in any common table of integrals. Nevertheless, we can find the imaginary part of $\log(g_{\alpha}(-\ir))$ by using the following lemma.
\begin{lemma} For $\mu\in\mathbb R$, let
\begin{equation}\label{iofmu}
J(\mu):=\int_0^{\infty}\log(1+t^{-2\mu})\frac{\er^{\ir \pi/2\mu}}{t^2-\er^{\ir \pi/\mu}}\,\dr t.
\end{equation}
Then for each $\mu>1/2$, which corresponds exactly to $\alpha<\pi/2$,
\begin{equation}
\Re(J(\mu))=\frac{\pi^2}{4}(1-\mu).
\end{equation}
\end{lemma}
\begin{proof} Change variables in $J(\mu)$ by letting $t=r^{-1}$. Then $dt=-r^{-2}dr$, and we have
\begin{equation}
J(\mu)=\int_0^{\infty}\log(1+r^{2\mu})\frac{\er^{\ir \pi/2\mu}}{1-r^2\er^{\ir \pi/\mu}}\,\dr r=-\int_0^{\infty}\log(1+r^{2\mu})\frac{\er^{-\ir \pi/2\mu}}{r^2-\er^{-\ir \pi/\mu}}\,\dr r.
\end{equation}
Breaking up the logarithm by bringing out $2\mu\log r$ gives
\begin{equation}
J(\mu)=-\int_0^{\infty}\log(1+r^{-2\mu})\frac{\er^{-\ir \pi/2\mu}}{r^2-\er^{-\ir \pi/\mu}}\,\dr r-2\mu\int_0^{\infty}\log(r)\frac{\er^{-\ir \pi/2\mu}}{r^2-\er^{-\ir \pi/\mu}}\,\dr r.
\end{equation}
But we now recognise the first integral as the complex conjugate of $J(\mu)$. Therefore
\begin{equation}\label{eq:simpleint}
\Re(J(\mu))=\frac{1}{2}(J(\mu)+\overline{J(\mu)})=-\mu\int_0^{\infty}\log(r)\frac{\er^{-\ir \pi/2\mu}}{r^2-\er^{-\ir \pi/\mu}}\,\dr r.
\end{equation}
The integral \eqref{eq:simpleint} may now be evaluated explicitly. It is the integral along the positive real axis of the function $\log(z)\dfrac{\er^{-\ir \pi/2\mu}}{z^2-\er^{-\ir \pi/\mu}}$, which has two poles (at $\pm \er^{-\ir \pi/2\mu}$) and a branch cut at $z=0$ which we take along the negative real axis. The decay at infinity and logarithmic growth at the origin enables us to use Cauchy's theorem to move the integral to an integral along the ray $\theta=(\pi/2-\pi/2\mu)$; since $\mu>1/2$, there are no singularities in the intervening region. Along this ray $z=\ir\rho \er^{-\ir \pi/2\mu}$, so we have
\begin{equation}
\Re(J(\mu))=\mu\int_0^{\infty}(\log\rho+\ir\frac{\pi}{2}(1-1/\mu))\frac{\ir}{\rho^2+1}\,\dr\rho.
\end{equation}
The imaginary part of this integral is, fortunately, zero, and the real part is $-\pi^2(\mu-1)/4$, completing the proof.
\end{proof}
As an immediate consequence of the Lemma, we have that for some function $E(\alpha)$,
\begin{equation}
g_{\alpha}(-\ir)=\er^{E(\alpha)}\er^{\ir (1-\mu)\pi/4}.
\end{equation}
This completes the proof of Theorem \ref{Peters}. \end{proof}

\subsection{Asymptotics of Peters solution as $z\to 0$} This analysis is substantially easier. First we have a small proposition:
\begin{proposition}\label{prop:asymptoticsofg} $g_{\alpha}(\zeta)\to 1$ as $|\zeta|\to\infty$, and the convergence is uniform in any sector away from the negative real axis.
\end{proposition}
\begin{proof} This is claimed on p. 329 of \cite{peters} and is immediate from \eqref{eq:g1} when $\Re(\zeta)>0$, but is not so clear for other values of $\arg(\zeta)$, so we rewrite the proof here.

From the definition of $g_{\alpha}(\zeta)$ it suffices to show that $I_{\alpha}(\zeta)\to 0$ as $|\zeta|\to\infty$. This is immediately clear when $\arg(\zeta)\in(-\alpha,\alpha)$, from the definition of $I_{\alpha}(\zeta)$ and a change of variables. For other values of $\arg(\zeta)$ we can use other representations of $I_{\alpha}(\zeta)$. For example, the representation \eqref{largerarearep} is good beyond $\zeta=\pm\pi/2$, and each term approaches zero as $|\zeta|\to\infty$ as long as $\zeta \er^{-\ir \alpha}$ is not on $M_2$ or $L_0$. The rate of convergence depends only on the distance of $\zeta/|\zeta|$ from $M_2$ and $L_0$. Continuing to move the contours as in the appendix of \cite{peters}, we can show that $I_{\alpha}(\zeta)$ goes to zero as $|\zeta|\to\infty$ as long as $\zeta$ is not on the branch cut. This completes the proof.
\end{proof}
Now recall that
\[f(z)=\ir\int_P\frac{g_{\alpha}(\zeta)}{\zeta+\ir}\er^{z\zeta}\,\dr\zeta.\]
Assume $|z|<1$. Change variables to $w=z\zeta$, then deform the contour back to $\arg(z)P$ (this works since all singularities of the integrand are inside $P$). We get
\[f(z)=\ir\int_{\arg(z)P}\frac{g_{\alpha}(w/z)}{w+\ir z}\er^w\,\dr w.\]
As $|z|\to 0$, $g_{\alpha}(w/z)$ converges uniformly to 1, and $(w+\ir z)^{-1}$ converges uniformly to $w^{-1}$. So as $|z|\to 0$,
\[f(z)\to \ir\int_{\arg(z)P}w^{-1}\er^w\,\dr w=-2\pi.\]
Thus we have established that $\lim_{z\to 0}f(z)$ exists and equals $-2\pi$. 

When a constant coefficient differential operator $\mathcal P$ of order $k$ is applied to $f(z)$, the analysis may be done similarly. The action of $\mathcal P$ brings down a factor of $\zeta^k$ inside the integral, which becomes $(w/z)^{k}$. Pulling the $z^{-k}$ out of the integral creates $\rho^{-k}$, and the remaining integral approaches the integral of $w^{k-1}\er^w$, which is zero. This completes 
the proof of Theorem \ref{wedgemodels} in the Robin-Neumann case.

\subsection{Robin-Dirichlet problem} This is also in \cite{peters}, again with some details but not others. We again look for a function $f(z)$ analytic in $\Sct{\alpha}$. The condition at $z\in \Sct{\alpha}\cap\{\theta=-\alpha\}$ is now $\Re(f(z))=0$. Following Peters, using the same notation for $g_{\alpha}(\zeta)$ as before, a solution is
\begin{equation}\label{eq:RDsolution}
f(z)=\ir\int_P\frac{g_{\alpha}(\zeta)}{\zeta+\ir}\zeta^{-\mu}\er^{z\zeta}d\zeta,\ k\in\mathbb Z.
\end{equation}
The only difference between the solutions is the multiplication by $\zeta^{-\mu}$ inside the integral.

We verify that \eqref{eq:RDsolution} satisfies the properties we need. Indeed it is obviously analytic in $\Sct{\alpha}$ (at least away from the corner point). The boundary condition on the real axis is verified in precisely the same way as for the Robin-Neumann problem, since $\zeta^{-\mu}$ is real when $\zeta$ is real. For the Dirichlet boundary condition, we again change variables to $w=\zeta \er^{-\ir \alpha}$, when \eqref{eq:RDsolution} becomes
\begin{equation}\label{eq:RDsolution1}
f(z)=\int_{P_1}\frac{w^{-\mu}\er^{\ir \alpha}g_{\alpha}(w\er^{\ir \alpha})}{w\er^{\ir \alpha}+\ir}\er^{pw}\,\dr w=\int_{P_1}w^{-\mu-1}\exp(-I_{\alpha}(w))\er^{\rho w}\,\dr w,
\end{equation}
where $P_1$ is symmetric with respect to the real axis and the branch cut is along the negative real axis. As for \eqref{eq:diff3}, this has zero real part, showing that \eqref{eq:RDsolution1} is a solution.

Now we analyse the asymptotics. Let
\[
\chi_{D,\alpha}=\chi+\frac{\pi^2}{4\alpha}=\frac{\pi}{4}\left(1+\frac{\pi}{2\alpha}\right).
\]
Analogously to the Robin-Neumann case, we have:
\begin{theorem}\label{PetersDirichlet} With $E(\alpha):(0,\pi/2)\to\mathbb R$ as in Theorem \ref{Peters}, there exists $R'_{\alpha}(z):W\times(0,\pi/2)\to\mathbb C$ such that
\begin{equation}
f(z)=\er^{E(\alpha)}\er^{-\ir (z-\chi_{D,\alpha})}+R'_{\alpha}(z),
\end{equation}
where for any fixed $\alpha\in(0,\pi/2)$, there exists a constant $C>0$ such that for all $z\in \Sct{\alpha}$,
\begin{equation}
|\Re(R'_{\alpha}(z))|\leq C|z|^{-2\mu},\ |\nabla_z (\Re(R'_{\alpha}(z))|\leq C|z|^{-2\mu-1}.
\end{equation}
\end{theorem}
\begin{proof} We mimic the proof of Theorem \ref{Peters}. By the same contour-deformation process, we have
\begin{equation}\label{eq:solution2mark1}
\begin{split}
f(z)&=2\mu^{1/2}g_{\alpha}(-\ir)\er^{-\ir z}\er^{\ir \pi\mu/2}+\frac{\mu^{1/2}}{\ir\pi}\int_{P'}\frac{g_{\alpha}(\zeta)\zeta^{-\mu}}{\zeta+\ir}\er^{z\zeta}\,\dr \zeta\\
&+\sum_{j=1}^m\er^{z\lambda_m}2\mu^{1/2}\mathop{\operatorname{Res}}\limits_{\zeta=\lambda_m}(g_{\alpha}(\zeta)\zeta^{-\mu}/(\zeta+\ir)).
\end{split}
\end{equation}
From the evaluation of $g_{\alpha}(-\ir)$ in the proof of Theorem \ref{Peters}, we immediately see that the first term is $\er^{E(\alpha)}\er^{-\ir (z-\chi_{D,\alpha})}$. The sum of residues in \eqref{eq:solution2mark1} has the needed decay estimates. And if we call the second term $\widetilde R'_{\alpha}(z)$, the methods of the proof of Theorem \ref{Peters} give
\begin{equation}\label{firstterm1b}
\widetilde R'_{\alpha}(z)=-\frac{\mu^{1/2}}{\pi}\int_{P'}\er^{z\zeta}\er^{-h(\zeta)}\frac{\dr\zeta}{\zeta},
\end{equation}
which has an extra factor of $\zeta^{-\mu}$ compared with \eqref{firstterm1a}. Changing variables as before gives
\begin{equation}\label{firstterm2b}
\widetilde R'_{\alpha}(z)=-\frac{\mu^{1/2}}{\pi}\int_{\arg(z)P'}\er^w\er^{-h(w/z)}\frac{\dr w}{w}.
\end{equation}
Plugging in $h(\cdot)=p(\cdot)+R(\cdot)$ gives
\begin{equation}\label{firstterm3b}
\widetilde R'_{\alpha}(z)=-\frac{\mu^{1/2}}{\pi}\int_{\arg(z)P'}\er^w\er^{-p(w/z)}\er^{-R(w/z)}\frac{\dr w}{w}.
\end{equation}
Using the estimates on $R(w/z)$, we have

\begin{equation}\label{firstterm4b}
\widetilde R'_{\alpha}(z)=-\frac{\mu^{1/2}}{\pi}\int_{\arg(z)P'}\er^w\er^{-p(w/z)}(1+O(|w/z|^{2\mu}))\frac{\dr w}{w}.
\end{equation}
For the part of the integral corresponding to the first term $1$, $p(w/z)$ is holomorphic in the full disk with no branch cut, so the contour may be deformed to a small disk around $w=0$ (plus a portion around any branch cut along the negative real axis starting from $w=-1$), then evaluated using residues. We see that the first term is equal to
\begin{equation}\label{firstterm2a}
-2\ir\mu^{1/2}+O(\er^{-|z|/C}),
\end{equation}
and therefore its real part is $O(\er^{-|z|/C})$, even better than needed. On the other hand, the exact same analysis as in the Robin-Neumann case shows that the second term is $O(|z|^{-2\mu})$, which is what we wanted. 

The gradient estimate follows as before, since differentiating \eqref{firstterm1b} just brings down an extra power of $\zeta$, which becomes $w/z$. The estimates on higher order derivatives follow as well.
\end{proof}

And the analysis of $f(z)$ as $z\to 0$ follows as before. We can show that
\[f(z)\sim z^{\mu}\int_{\arg(z) P}w^{-\mu-1}\er^w \dr w.\]
Since $\mu>0$, we have $\lim_{z\to 0}f(z)=0$, and moreover $f(z)=O(|z|^{\mu})$. As before, the action of any differential operator $\mathcal P$ brings down an extra factor of $\zeta^k=w^k/z^k$ and the analysis proceeds similarly. This completes the proof of Theorem \ref{wedgemodels} in the Robin-Dirichlet case.

\section{Proof of Proposition \ref{prop:odeasymp}.}
\label{naimarkproof}
\subsection{Plan of the argument}
Here we prove Proposition \ref{prop:odeasymp} by using the analysis of \cite{naimark}. As we will see, it is a special case of Theorem 2 in \cite{naimark}. In what follows we assume that 
$q$ is even, see Remark \ref{naimark:proof}. Throughout we let $n=2q$ and note that $n=0$ mod 4.

We first claim the following technical lemma:
\begin{lemma}\label{naimarklemma} Suppose $q\in\mathbb N$ is even. Then the boundary conditions in \eqref{odeneumann} and \eqref{odedirichlet} are regular in the sense of \cite[section 4.8]{naimark}. Moreover, using the notation of \cite[equations (41)--(42)]{naimark}, $\theta_0=0$, and $\theta_{1}=\theta_{-1}\neq 0$ in both cases.
\end{lemma}

Assuming this lemma, we obtain asymptotics for the eigenvalues. Indeed, let $\xi'$ and $\xi''$ be the roots of the equation $\theta_1\xi^2+\theta_0\xi+\theta_{-1}=0$; by our calcuation, these are the roots of $\xi^2+1=0$:
\[\xi'=-\ir,\xi''=\ir.\]
From \cite[Theorem 2, equations (45a) and (45b)]{naimark}, using the fact that $n=0$ mod 4, all sufficiently large eigenvalues form two sequences:
\[
(\lambda_k')^{2q}=(2k\pi)^{2q}\left(1-\frac{q\ln_0\xi'}{k\pi\ir}+O\left(\frac{1}{k^2}\right)\right);
\]
\[
(\lambda_k'')^{2q}=(2k\pi)^{2q}\left(1-\frac{q\ln_0\xi''}{k\pi\ir}+O\left(\frac{1}{k^2}\right)\right).
\]
Taking $2q$-th roots and using Taylor series, and plugging in $\xi'$ and $\xi''$:
\[
\lambda_k'=2k\pi\left(1-\frac{\ln_0\xi'}{2k\pi\ir}+O\left(\frac{1}{k^2}\right)\right)=2k\pi+\frac{\pi}{2}+O\left(\frac 1k\right);\]
\[\lambda_k''=2k\pi\left(1-\frac{\ln_0\xi''}{2k\pi\ir}+O\left(\frac{1}{k^2}\right)\right)=2k\pi-\frac{\pi}{2}+O\left(\frac 1k\right).\]
Combining the two sequences yields Proposition \ref{prop:odeasymp}.

\begin{remark}Theorem 2 in \cite{naimark}, in some editions and English translations, only states that there are sequences of eigenvalues of the desired form, rather than stating that \emph{all} sufficiently large eigenvalues form those sequences. However, the original Russian edition claims the stronger version, and it follows immediately from the proof in \cite{naimark} anyway.
\end{remark}

\subsection{Proof of Lemma \ref{naimarklemma}}

To find $\theta_1$, $\theta_0$, and $\theta_{-1}$, we need to analyse the notation of \cite{naimark}. The boundary conditions will be seen to be regular conditions of Sturm type, hence falling under the analysis of \cite[section 4.8]{naimark}, in particular giving $\theta_0=0$ and $\theta_1$ and $\theta_{-1}$ defined by (41)--(42) in \cite{naimark}. To analyse these determinants we need to identify the numbers $\omega_j$, $k_j$, and $k_j'$. We do this first. 

To analyse the $\omega_j$, see the discussion before \cite[(38)]{naimark}. Note that the sector $S_0$ referred to in the statement of theorem 2 is the sector of the complex plane with argument between $0$ and $\frac{\pi}{n}$ (see \cite[(3)]{naimark}). We then let $\omega_1,\dots,\omega_n$ be the $n$th roots of $-1$, arranged so that for $\rho\in S_0$,
\[\Re(\rho\omega_1)\leq\Re(\rho\omega_2)\leq\dots\leq\Re(\rho\omega_n).\]
To figure out this ordering, let
\[\omega:=\er^{\ir \pi/n}.\]
 Then, in order, the set $\{\omega_1,\dots,\omega_n\}$ is
\[\{\omega^{2q-1},\omega^{-(2q-1)},\dots,\omega^{3},\omega^{-3},\omega^1,\omega^{-1}\}.\]
Note further that $\omega_{\mu}=\omega_{q}$; since $q$ is even,
\[\omega_{\mu}=\omega^{-(q+1)}; \qquad \omega_{\mu+1}=\omega^{q-1}.\]
In particular we actually have $\omega_{\mu+1}=-\omega_{\mu}$. This will be useful.

Now we need to understand the $k_j$ and $k_j'$, which is done via a direct comparison with \cite[(40)]{naimark}. The sums all vanish, and only the first terms remain. We see that we always have $k_j'=k_j$ in both the Neumann and Dirichlet settings, and for Neumann we have
\[k_1=2q-1,\ k_2=2q-2,\dots,k_q=q.\]
For Dirichlet we have
\[k_1=q-1,\ k_2=q-2,\dots,k_q=0.\]
So these conditions are of Sturm type.

Since the conditions are Sturm type we immediately have $\theta_0=0$ as well as the formulas \cite[(41)--(42)]{naimark} for $\theta_{-1}$ and $\theta_1$. Assume for the moment we are working in the Neumann setting. Define the matrices
\[
A_q:=\begin{bmatrix}
\omega_1^{k_1} & \cdots & \omega_{\mu}^{k_1}\\ 
\vdots & \ddots & \vdots \\ 
\omega_1^{k_{\mu}} & \cdots & \omega_{\mu}^{k_{\mu}}
\end{bmatrix}
=
\begin{bmatrix}
\omega^{(2q-1)(2q-1)} & \omega^{-(2q-1)(2q-1)} & \cdots & \omega^{-(q+1)(2q-1)}\\ 
\vdots & \vdots & \cdots & \vdots \\ 
\omega^{(2q-1)(q)} & \omega^{-(2q-1)(q)} & \cdots & \omega^{-(q+1)(q)}
\end{bmatrix},
\]
\[
B_q:=\begin{bmatrix}
\omega_{\mu +1}^{k_1} & \cdots & \omega_{n}^{k_1}\\ 
\vdots & \ddots & \vdots \\ 
\omega_{\mu +1}^{k_{\mu}} & \cdots & \omega_{n}^{k_{\mu}}
\end{bmatrix}
=
\begin{bmatrix}
\omega^{(q-1)(2q-1)} & \omega^{-(q-1)(2q-1)} & \cdots & \omega^{-1(2q-1)}\\ 
\vdots & \vdots & \cdots & \vdots \\ 
 \omega^{(q-1)(q)} & \omega^{-(q-1)(q)} & \cdots & \omega^{-1(q)}
 \end{bmatrix}.
 \]
Let $A_q'$ be $A_q$ with the last column replaced by the first column of $B_q$, and let $B_q'$ be $B_q$ with the first column replaced by the last column of $A_q$. Note also that since $\omega_{\mu+1}=-\omega_{\mu}$, those switched columns are the same except that the entry in each odd row is multiplied by $-1$. In particular their last entries are the same, since $q$ is even. Then \cite[(41)--(42)]{naimark} give
\begin{equation}\label{thetas}
\theta_{-1}=\pm\det A_q\det B_q;\ \theta_1=\pm\det A_q'\det B_q',
\end{equation}
where the signs are the same (they both come from the same equation before (41) and (42) in \cite{naimark}). Thus it suffices to show that
\begin{equation}\label{determinantident}\frac{\det A_q'}{\det A_q}=\frac{\det B_q'}{\det B_q}.\end{equation}

The trick will be to reduce these to Vandermonde determinants by taking a factor out of each column; the factor we take out will always be the \emph{last} entry in each column, which is simplified by the fact that $w^q=\ir$. For example,
\[\det B_q=\ir^{(q-1)}\ir^{-(q-1)}\ir^{(q-3)}\ir^{-(q-3)}\dots \ir^1\ir^{-1}\det D_q=\det D_q,\]
where
\[D_q=\begin{bmatrix}\omega^{(q-1)(q-1)} & \omega^{-(q-1)(q-1)} & \cdots &\omega^{-1(q-1)}\\ \vdots & \vdots & \cdots & \vdots \\ \omega^{(q-1)(1)} & \omega^{-(q-1)(1)} & \cdots & \omega^{-1(1)}\\ 1 & 1 & \cdots & 1\end{bmatrix}.\]
Similarly, $\det A_q=\det C_q$, where $C_q$ is obtained as with $D_q$. As for $B_q'$ and $A_q'$, note that the last entry of the last column of $A_q$ is the same as the last entry of the first column of $B_q$, so all the pre-factors still cancel and we have $\det A_q'=\det C_q'$, $\det B_q'=\det D_q'$, with $D_q'$ and $C_q'$ obtained by switching the first column of $D_q$ for the last column of $C_q$.

Switching the order of all of the rows introduces a factor of $(-1)^{q/2}$ in each determinant and brings us to a Vandermonde matrix in each case:
\begin{align*}
\det B_q&=(-1)^{q/2}\det V(\omega^{q-1},\omega^{-(q-1)},\dots,\omega^1,\omega^{-1});\\
\det B_q'&=(-1)^{q/2}\det V(\omega^{-(q+1)},\omega^{-(q-1)},\dots,\omega^1,\omega^{-1});\\
\det A_q&=(-1)^{q/2}\det V(\omega^{2q-1},\omega^{-(2q-1)},\dots,\omega^{q+1},\omega^{-(q+1)});\\
\det A_q'&=(-1)^{q/2}\det V(\omega^{2q-1},\omega^{-(2q-1)},\dots,\omega^{q+1},\omega^{q-1}).
\end{align*}
Each of these Vandermonde determinants can be computed explicitly via the well-known pairwise difference formula. Considering $\det B_q'/\det B_q$, all the pairwise differences cancel except for those involving the first element. The first element just changes by an overall sign from $B_q$ to $B_q'$, and we get
\[\frac{\det B_q'}{\det B_q}=\frac{(\omega^{-1}+\omega^{q-1})(\omega^{1}+\omega^{q-1})\dots(\omega^{-(q-1)}+\omega^{q-1})}{(\omega^{-1}-\omega^{q-1})(\omega^{1}-\omega^{q-1})\dots(\omega^{-(q-1)}-\omega^{q-1})}.\]
Similarly,
\[\frac{\det A_q'}{\det A_q}=\frac{(-\omega^{q-1}-\omega^{q+1})\dots(-\omega^{q-1}-\omega^{-(2q-1)})(-\omega^{q-1}-\omega^{2q-1})}{(\omega^{q-1}-\omega^{q+1})\dots(\omega^{q-1}-\omega^{-(2q-1)})(\omega^{q-1}-\omega^{2q-1})}.\]
Pulling a minus sign out of each term on the top and the bottom and reversing the order of the multiplication gives
\[\frac{\det A_q'}{\det A_q}=\frac{(\omega^{q-1}+\omega^{2q-1})(\omega^{q-1}+\omega^{-(2q-1)})\dots(\omega^{q-1}+\omega^{q+1})}{(-\omega^{q-1}+\omega^{2q-1})(-\omega^{q-1}+\omega^{-(2q-1)})\dots (-\omega^{q-1}+\omega^{q+1})}.\]
Flipping the order of addition within each term:
\[\frac{\det A_q'}{\det A_q}=\frac{(\omega^{2q-1}+\omega^{q-1})(\omega^{-(2q-1)}+\omega^{q-1})\dots(\omega^{q+1}+\omega^{q-1})}{(\omega^{2q-1}-\omega^{q-1})(\omega^{-(2q-1)}-\omega^{q-1})\dots (\omega^{q+1}-\omega^{q-1})}.\]
Using the fact that $\omega^{2q}=1$ to simplify one element in each term, we see that this fraction is precisely $\frac{\det B_q'}{\det B_q}$, completing the proof of Lemma \ref{naimarklemma} for the Neumann case.

The Dirichlet case is very similar, in fact easier; we have $k_1$ through $k_q$ each decreasing by $q$. In particular the last entries of each column are $1$, so the matrices corresponding to $A_q, B_q, A_q',B_q'$ are already Vandermonde. In fact, they are exactly equal to $C_q, D_q, C_q'$, and $D_q'$ respectively, so our previous computation completes the proof in this case as well.
\qed

\section{Numerical examples}\label{app:examples} In this Appendix we present some numerical results. 
All computations have been performed using the package FreeFem++ \cite{FreeFem}.

\subsection{Example illustrating Theorems \ref{thm:main} and \ref{thm:Dirichlet}}\label{sec:ex1}

Let $\Omega=\bigtriangleup ABZ$ be a triangle with $L=1$, $\alpha=\frac{2 \pi}{5}$ and $\beta=\frac{\pi}{6}$. We consider sloshing with Neumann or Dirichlet conditions on $\mathcal{W}=[A,Z]\cup[Z,B]$. The eigenvalues and quasi-frequencies are given in the table below, and we see that the error is indeed quite small in both Neumann and Dirichlet cases:

\begin{figure}[htb!]
\begin{center}
\includegraphics[width=0.7\textwidth]{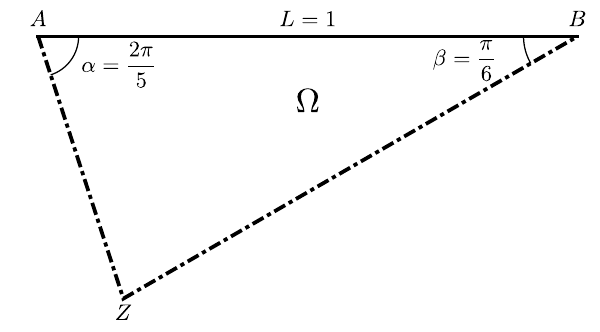}
\end{center}
\caption{Triangle from Example \ref{sec:ex1}\label{fig:example1}}
\end{figure}

\begin{center}
\begin{tabular}{@{}>{$}r<{$}*{6}{>{$}c<{$}}@{}}
\toprule
k&\lambda_{k}&\sigma_{k}&\left|\dfrac{\sigma_{k}}{\lambda_{k}}-1\right|&\lambda_{k}^D&\sigma_{k}^D&\left|\dfrac{\sigma_{k}^D}{\lambda_{k}^D}-1\right|\\[1ex]\midrule
1 & 0. & -0.88357 & \ & 2.43592 & 2.45437 & 7.58\times 10^{-3} \\\hline
 2 & 0.85626 & 0.68722 & 1.97\times 10^{-1} & 4.02389 & 4.02517 & 3.17\times 10^{-4} \\\hline
 3 & 2.28840 & 2.2580 & 1.33\times 10^{-2} & 5.59623 & 5.59596 & 4.83\times 10^{-5} \\\hline
 4 & 3.82292 & 3.8288 & 1.54\times 10^{-3} & 7.16681 & 7.16676 & 7.11\times 10^{-6} \\\hline
 5 & 5.39779 & 5.3996 & 3.37\times 10^{-4} & 8.73757 & 8.73755 & 1.81\times 10^{-6} \\\hline
 6 & 6.96977 & 6.9704 & 9.22\times 10^{-5} & 10.3084 & 10.3084 & 1.44\times 10^{-6} \\\hline
 7 & 8.54086 & 8.5412 & 4.00\times 10^{-5} & 11.8792 & 11.8791 & 1.74\times 10^{-6} \\\hline
 8 & 10.1118 & 10.112 & 2.03\times 10^{-5} & 13.4500 & 13.4499 & 2.36\times 10^{-6} \\\hline
 9 & 11.6827 & 11.683 & 1.11\times 10^{-5} & 15.0208 & 15.0207 & 3.26\times 10^{-6} \\\hline
 10 & 13.2535 & 13.254 & 5.90\times 10^{-6} & 16.5916 & 16.5915 & 4.47\times 10^{-6} \\\bottomrule
\end{tabular}
\end{center}

\subsection{Example illustrating Remark \ref{rem:believe}}

We define two domains $\Omega_\pm$ by setting a \emph{curved} sloshing surface 
\[
S_\pm=\left\{\left(x,\pm\frac{1}{2\pi}\sin(2\pi x)\right)\mid 0<x<1\right\},
\]
so that the length of the  sloshing surface is
\[
L=\int_0^1 \sqrt{1+ \cos^2(2\pi x)}\,\dr x=\frac{2\sqrt{2}}{\pi}E\left(\frac12\right)\approx 1.21601,
\]
where $E$ denotes a complete elliptic integral of the second kind. Let
\[
\begin{split}
\W_1&=\W_{1,\pm}=\left\{z:\ |z-1/2|=1/2, -\pi\le\arg(z)\le-\pi/2\right\}, \\
\W_2&=\W_{2,\pm}=\left\{z:\ |z-1/2|=1/2, -\pi/2\le\arg(z)\le 0\right\},
\end{split}
\]
so that 
\[
\alpha_+=\beta_-=\frac{3\pi}{4}, \qquad \alpha_-=\beta_+=\frac{\pi}{4}.
\]
In both cases the Dirichlet boundary condition is imposed on $\W_1$ and the Neumann one on $\W_2$.

\begin{figure}[hbt!]
\begin{center}
\includegraphics{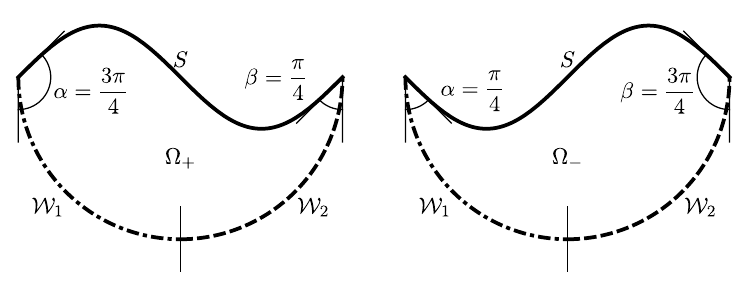}
\end{center}
\caption{Domains $\Omega_+$ (left) and  $\Omega_-$ (right)\label{fig:example2}}
\end{figure}

The quasifrequencies  are then given by \eqref{mixedasymp}, which after omitting $O-$terms simplifies to
\[
\sigma_{+,k}=\frac{\pi}{L}\left(k-\frac16\right),\qquad \sigma_{-,k}=\frac{\pi}{L}\left(k-\frac56\right).
\] 
The comparisons between the numerically calculated eigenvalues and the quasifrequencies are below. Note that there is very good agreement even for low $k$.

\begin{center}
\begin{tabular}{@{}>{$}r<{$}*{6}{>{$}c<{$}}@{}}
\toprule
k&\lambda_{+,k}&\sigma_{+,k}&\left|\dfrac{\sigma_{+,k}}{\lambda_{+,k}}-1\right|&\lambda_{-,k}&\sigma_{-,k}&\left|\dfrac{\sigma_{-,k}}{\lambda_{-,k}}-1\right|\\[1ex]\midrule
 1 & 1.02371 & 2.15294 & 1.10  & 1.24543 & 0.430589 & 6.54\times 10^{-1} \\\hline
 2 & 5.65749 & 4.73648 & 1.63\times 10^{-1} & 2.63524 & 3.01412 & 1.44\times 10^{-1} \\\hline
 3 & 8.13194 & 7.32001 & 9.98\times 10^{-2} & 5.55627 & 5.59765 & 7.45\times 10^{-3} \\\hline
 4 & 10.3085 & 9.90354 & 3.93\times 10^{-2} & 8.22122 & 8.18119 & 4.87\times 10^{-3} \\\hline
 5 & 12.8138 & 12.4871 & 2.55\times 10^{-2} & 10.6845 & 10.7647 & 7.51\times 10^{-3} \\\hline
 6 & 15.3856 & 15.0706 & 2.05\times 10^{-2} & 13.1122 & 13.3483 & 1.80\times 10^{-2} \\\hline
 7 & 17.9151 & 17.6541 & 1.46\times 10^{-2} & 15.7600 & 15.9318 & 1.09\times 10^{-2} \\\hline
 8 & 20.4310 & 20.2377 & 9.46\times 10^{-3} & 18.4111 & 18.5153 & 5.66\times 10^{-3} \\\hline
 9 & 22.9800 & 22.8212 & 6.91\times 10^{-3} & 21.0011 & 21.0988 & 4.66\times 10^{-3} \\\hline
 10 & 25.5511 & 25.4047 & 5.73\times 10^{-3} & 23.5873 & 23.6824 & 4.03\times 10^{-3} \\\bottomrule
\end{tabular}
\end{center}

\end{document}